\documentclass[11pt]{article}
\usepackage{amsmath,amssymb,amsthm,graphicx,pst-tree}
\usepackage[sort,round]{natbib}

\title{Algebraic Factor Analysis: Tetrads, Pentads and Beyond}
\author{Mathias Drton, Bernd Sturmfels and Seth Sullivant}
\date{\today}

\DeclareMathOperator{\codim}{codim}

\DeclareMathOperator{\rank}{rank}
\DeclareMathOperator{\trace}{tr}
\DeclareMathOperator{\E}{E}
\DeclareMathOperator{\Var}{Var}

\newcommand{\RRR}{\mathbb{R}}
\newcommand{\CCC}{\mathbb{C}}

\newcommand{\ND}{\mathcal{N}}
\newtheorem{theorem}{Theorem}
\newtheorem{corollary}[theorem]{Corollary}
\newtheorem{proposition}[theorem]{Proposition}

\newtheorem{conjecture}[theorem]{Conjecture}
\theoremstyle{definition}

\newtheorem{remark}[theorem]{Remark}
\newtheorem{example}[theorem]{Example}
\newtheorem{question}[theorem]{Question}

\begin{document}
\maketitle


\begin{abstract}
\noindent
Factor analysis refers to a statistical model in which observed variables
are conditionally independent given fewer hidden variables, known as
factors, and all the random variables follow a multivariate normal
distribution.  The parameter space of a factor analysis model is a subset
of the cone of positive definite matrices. This parameter space is studied
from the perspective of computational algebraic geometry.  Gr\"obner bases
and resultants are applied to compute the ideal of all polynomial functions
that vanish on the parameter space.  These polynomials, known as model
invariants, arise from rank conditions on a symmetric matrix under
elimination of the diagonal entries of the matrix.  Besides revealing the
geometry of the factor analysis model, the model invariants also furnish
useful statistics for testing goodness-of-fit.
\end{abstract}


\section{Introduction}
\label{sec:intro}

In factor analysis, correlated continuous variables are modeled as
conditionally independent given hidden (latent) variables that are termed
{\em factors}.  Sometimes factor analysis serves as a tool for
dimension-reduction; the possibly many observed variables are summarized by
fewer factors.  However, in many applications the focus is on interpreting
the factors as unobservable theorized concepts.  In fact, the desire to
explain observed correlations between individuals' exam performances by the
concept of intelligence was the driving force in the original development
of factor analysis \citep{spearman:04,spearman:27}.

Currently, statistical inference in factor analysis is often based
exclusively on parametric representations and on maximum likelihood
estimates computed in iterative procedures such as the EM algorithm
\citep{rubin:82}.  In the early days of factor analysis, however, much
attention was directed to {\em model invariants\/}, that is, to polynomial
equality relations that the model imposes on the entries of the covariance
matrix of the observed variables. We refer to \citet[Sect.~1]{harman:76}
for some history.  It should be noted that factor analysis also leads to
inequality constraints \citep[e.g.][p.~117]{bekker:87,harman:76}, which we
do not address here.  The best known invariants are the {\em tetrads\/},
also called {\em tetrad differences\/}, which arise in one-factor models.
The name tetrad reflects that the polynomial arises in one-factor analysis
with four observed variables.  For example, if
$\Psi=(\psi_{ij})\in\RRR^{4\times 4}$ is a covariance matrix in the
parameter space of a one-factor analysis model, then there are, up to sign
change, three tetrads, namely
\begin{equation}
  \label{eq:tetrad}
  \psi_{12}\psi_{34}-\psi_{14}\psi_{23}, \quad 
  \psi_{14}\psi_{23}-\psi_{13}\psi_{24}, \quad 
  \psi_{12}\psi_{34}-\psi_{13}\psi_{24},
\end{equation}
and all three tetrads evaluate to zero.  Tetrads have played a major role
throughout the history of factor analysis.  They also appear in recent
research, for example, in work on model identifiability \citep{grzebyk:04}
and on dichotomized Gaussian models for multivariate binary variables
\citep{cox:02}.  While tetrads are ubiquitous in the literature, there has
been very little work attempting to find invariants of models with more
than one factor.  The work by \cite{kelley:35} who derived the {\em
  pentad\/}, a fifth degree polynomial vanishing over covariance matrices
from two-factor models, constitutes the exception.  Since then virtually no
progress has been made towards determining the invariants of factor
analysis models.  \citet[p.~77]{harman:76} summarizes the state of
knowledge as follows: {\em ``When the number of factors is greater than two
  the work of computing determinants of the fourth or higher order becomes
  so laborious that no explicit conditions corresponding to the tetrads or
  pentad criterion have been worked out.''}

Computational difficulties aside, the apparent ease of data analysis using
solely the parametric model representation has inhibited progress on the
determination of higher-order invariants.  However, parametric approaches
are not without their problems.  On one hand, the likelihood function of a
factor analysis model may have multiple local maxima \citep{rubin:82},
rendering its maximization difficult.  On the other hand, the use of
information criteria such as BIC in exploratory factor analysis is
complicated by the presence of singularities \citep{geiger:01}.  We believe
that a better mathematical understanding of factor analysis models will be
helpful in addressing these issues.  One step in this direction is the work
of \citet{ellis:04} who applied algebraic topology to study singularities
arising in factor analysis.  Our interest lies in the algebraic geometry
that expresses itself in the model invariants.  This has a pragmatic side
because the invariants can serve as useful statistics for testing model fit
and for constraint-based model selection.  Both the desire to find new test
statistics as well as the wish for an understanding of the geometry of
factor analysis constitute the motivation for this paper.

The paper is outlined as follows. We begin with a review of the
factor analysis model (Section \ref{sec:factoranalysis}) and discuss
the use of invariants as test statistics (Section \ref{sec:tetradsstats}).
In Section \ref{sec:algsetup} we place the problem of determining
invariants in the framework of algebraic statistics \citep{ascb, pistone:01}.
In fact, this is one of the first studies in algebraic statistics which
deals with continuous rather than discrete random variables.

We shall see in Section \ref{sec:ideals} that the tetrads form a Gr\"obner
basis of the ideal of invariants of a one-factor model.  This implies that
all invariants of one-factor models can be written as polynomial
combinations of tetrads, which is claimed correctly but without proof in
\citet[p.~84]{glymour:87}.  For models with two or more factors, we
performed extensive Gr\"obner basis computations, using {\tt Macaulay 2}
and {\tt Singular}, for factor analysis models with up to nine observed
variables and up to five factors.  In Section \ref{sec:resultants}, we show
that multilinear resultants provide a useful method of finding individual
invariants even when the whole ideal of invariants cannot be determined.

Our computational experiments lead to a series of conjectures and problems
presented in Section \ref{sec:conjectures}.  In particular, we conjecture
for two-factor models that $3\times 3$-minors and pentads generate the
ideal of invariants. For models with arbitrarily many factors we conjecture
that the ideal is generated by polynomials arising from consideration of
submatrices whose size depends on the number of factors but not on the
number of observed variables.  We believe that these conjectures are of
independent interest for commutative algebra.  In Section \ref{sec:future}
we propose some future research directions of statistical interest.


\section{Factor analysis}
\label{sec:factoranalysis}

\begin{figure}[tb]
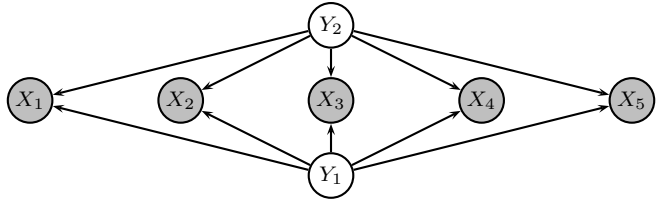

  \begin{minipage}{13cm}
    \vspace{5cm}\hspace{1cm}
    \scriptsize
    \psset{linewidth=0.8pt}          
    \newlength{\MyLength}
    \settowidth{\MyLength}{$X$}
    \newcommand{\myNode}[2]{\circlenode{#1}{\makebox[\MyLength]{#2}}}
    \psset{fillcolor=lightgray, fillstyle=solid}          
    \rput(2,3){\myNode{1}{$X_1$}}
    \rput(4,3){\myNode{2}{$X_2$}}
    \rput(6,3){\myNode{3}{$X_3$}}
    \rput(8,3){\myNode{4}{$X_4$}}
    \rput(10,3){\myNode{5}{$X_5$}}
    \psset{fillstyle=none}          
    \rput(6,2){\myNode{y1}{$Y_1$}}
    \ncline{->}{y1}{1}
    \ncline{->}{y1}{2}
    \ncline{->}{y1}{3}
    \ncline{->}{y1}{4}
    \ncline{->}{y1}{5}
    \psset{fillstyle=solid, fillcolor=white}          
    \rput(6,4){\myNode{y2}{$Y_2$}}
    \ncline{->}{y2}{1}
    \ncline{->}{y2}{2}
    \ncline{->}{y2}{3}
    \ncline{->}{y2}{4}
    \ncline{->}{y2}{5}
    \end{minipage}
    \vspace{-1cm}
  \caption{Graphical representation of the factor analysis model
    $\mathbf{F}_{5,2}$.} 
    \label{fig:factordag}
\end{figure}

Factor analysis concerns a Gaussian hidden variable model with $p$ observed
variables $X_i$, where $i \in [p] = \{1,\ldots,p\}$, and $m$ hidden
variables $Y_j$, where $j \in [m] = \{1,\ldots,m\}$. It is assumed that
$(X,Y)$ follows a joint multivariate normal distribution with positive
definite covariance matrix.  The {\em factor analysis model}
$\,\mathbf{F}_{p,m}$ is defined by the requirement that the observed
variables $X_i$, $i\in [p]$, are conditionally independent given the hidden
variables $Y_j$, $j\in [m]$.  The factor analysis model $\mathbf{F}_{p,m}$
can be visualized using the graphical model formalism \citep{lau:bk}, in
which the dependence structure between observed and hidden variables is
encoded by an acyclic directed graph.  This directed graph has the vertex
set $\{X_1,\ldots,X_p, Y_1, \ldots,Y_m\}$, and the edges are $Y_j\to X_i$
for all $j\in [m]$ and $i\in [p]$, as shown in Figure \ref{fig:factordag}
for $m = 2$ and $p = 5$.

We start out by deriving the following parametric representation of our
model.

\begin{proposition} \label{prop:factordef}
  The factor analysis model $\mathbf{F}_{p,m}$ is the family of
  multivariate normal distributions $\ND_p(\mu,\Psi)$ on $\RRR^p$ whose
  mean vector $\mu$ is an arbitrary vector in $\RRR^p$ and whose
  covariance matrix $\Psi$ lies in the (non-convex) cone
  \begin{equation}
    \label{eq:F}
    \begin{split}
      F_{p,m} \,\, & = \,\,\{ \Sigma+\Lambda\Lambda^t \in\RRR^{p\times
        p}\,:\, \Sigma>0 \text{ diagonal}, \; \Lambda\in \RRR^{p\times m}
        \}\\ 
      &= \,\,\{ \,\Sigma+\Gamma\in\RRR^{p\times
        p}\,\,:\, \Sigma>0
      \text{ diagonal},\; 
      \Gamma\ge 0\text{ symmetric},\; \text{rank}(\Gamma)\le m \}.
    \end{split}  
  \end{equation}
\end{proposition}
Here the notation $A >0$ means that $A$ is a positive definite matrix
(i.e.,~all eigenvalues are positive), and similarly $ A \ge 0$ means that
$A$ is a positive semidefinite matrix.

\begin{proof}
Consider the joint covariance matrix of the
fully observed model underlying $\mathbf{F}_{p,m}$,
\begin{equation}
\label{jointcovariance}
 \mathrm{Cov}\begin{pmatrix} X\\ Y\end{pmatrix} \quad = \quad   
\begin{pmatrix}
\Psi & \Lambda \\
\Lambda^t & \Phi
\end{pmatrix}.
\end{equation}
The entries of this matrix are constrained by the
conditional independence statements
\begin{equation}
\qquad \qquad
X_i \perp \!\!\! \perp X_j \, | \, \{Y_1,Y_2,\ldots,Y_m \} 
\qquad \qquad ( 1 \leq i < j \leq p),
\end{equation}
which translate into the vanishing of
the corresponding $(m+1) \times (m+1)$-determinants: 
\begin{equation}
\label{mplusonedet}  \qquad \qquad
{\rm det} \begin{pmatrix}
 \psi_{ij}       & \Lambda_{i*}   \\
 \Lambda_{j*}^t    & \Phi   
 \end{pmatrix} \quad = \quad  
 \psi_{ij} \cdot {\rm det}(\Phi) \,-\,
 \Lambda_{i *} \cdot {\rm adj}(\Phi) \cdot \Lambda_{j*}^t
  \quad = \quad 0.
\end{equation}
We refer to \citet{matus:05} for a general discussion
on how to translate
conditional independence statements for
Gaussian random variables into polynomial algebra.

The determinantal constraint (\ref{mplusonedet}) allows us to
block-diagonalize the positive definite matrix (\ref{jointcovariance}) as
follows:
\begin{equation}
\begin{pmatrix}
\Sigma & 0 \\
0 & \Phi 
\end{pmatrix}
\quad  = \quad
\begin{pmatrix}
    I_p & - \Lambda \Phi^{-1} \\
    0 & I_m
\end{pmatrix}
\cdot
\begin{pmatrix}
\Psi & \Lambda \\
\Lambda^t & \Phi
\end{pmatrix}
\cdot
\begin{pmatrix}
    I_p & 0\\
    - \Phi^{-1} \Lambda^t & I_m
\end{pmatrix}.
\end{equation}
Upon multiplication by ${\rm det}(\Phi) > 0$,
the entry of the matrix  $\,\Sigma \, = \, \Psi - 
 \Lambda \cdot \Phi^{-1} \cdot \Lambda^t\,$
 in row $i$ and column $j$ is equal to (\ref{mplusonedet}),
 so this positive definite matrix is diagonal if and only if
 $X$ satisfies the model $\mathbf{F}_{p,m}$.
 This holds if and only if its
 covariance matrix $\Psi$ has the form $\,\Psi \,=\,\Sigma+ 
 \Lambda \cdot \Phi^{-1} \cdot \Lambda^t \, $
 if and only if $\Psi$ is in the cone $F_{p,m}$. 
 \end{proof}

 In what follows we generally identify the factor analysis model
 $\mathbf{F}_{p,m}$ with its parameter space $F_{p,m}$.  The description
 given in Proposition \ref{prop:factordef} shows that $F_{p,m}$ is a
 parametrically presented subset of the space $ \RRR^{\binom{p+1}{2}}$ of
 symmetric $p \times p$-matrices.  The {\em dimension} $d = \dim(F_{p,m})$
 of the model $F_{p,m}$ is the maximal rank of the Jacobian matrix of that
 parametrization. The {\em codimension} of $F_{p,m}$ is $\binom{p+1}{2} -
 d$.

\begin{theorem}
  \label{thm:dim}
  The dimension and the codimension of the factor analysis model are
  \[
  \dim(F_{p,m})=\min\left\{p(m+1)-
    \binom{m}{2},\binom{p+1}{2}\right\},
  \]  
 \[
  \codim(F_{p,m})=\max\left\{\binom{p-m}{2} -m,0\right\}.
  \]
  Thus the codimension of the factor analysis model is positive if and only if
  \begin{equation}
  \label{ineqSolveForp}
  p \, \ge \, \left\lfloor m+ \frac{1}{2}\sqrt{8m+1}+\frac{1}{2}
  \right\rfloor + 1.
  \end{equation}
\end{theorem}

\begin{proof}
  Using orthogonal transformations as in the
  QR-decomposition, every $\Psi\in F_{p,m}$ can be written as
  $\Psi = \Sigma+\Lambda\Lambda^t$ with $\Lambda=(\lambda_{ij})$ being
  lower-triangular in the sense that
  \[
  \Lambda \,\, \in \,\, L_{p,m} \,=\, \left\{ \Lambda\in\RRR^{p\times m}\mid
    \lambda_{ij} = 0 \;\text{for all}\; 1\le i<j\le m \right\};
  \]
  see also \citet[pp.~121 and 124]{anderson:56}.  
  Thus the factor analysis model $F_{p,m}$ is the image
  of the following polynomial map:
  \begin{equation}
  \label{reducedpara}
  \RRR^m_{> 0} \times L_{p,m} \,\rightarrow \, \RRR^{\binom{p+1}{2}} \,, \,\,\,
  (\Sigma,\Lambda) \,\mapsto\, \Sigma+\Lambda\Lambda^t.
  \end{equation}
  The coordinates of the parametrization (\ref{reducedpara}) are
  \[
  \psi_{ij} \,\,\, = \,\,\, \begin{cases} \,
    \sigma_{ii}+\sum_{r=1}^{\min(i,m)} \lambda_{ir}^2 &\text{if}\; i=j,\\
 \,   \sum_{r=1}^{\min(i,m)} \lambda_{ir}\lambda_{jr} &\text{if}\;
    i <j.
    \end{cases}
  \]
  The dimension of the domain and the image space of
  (\ref{reducedpara}) are $\,p(m+1) - \binom{m}{2}\,$ and
  $\,\binom{p+1}{2}\,$ respectively, so the minimum of these two
  numbers is an upper bound for $F_{p,m}$.  To prove that this upper
  bound is tight, we will show that the Jacobian matrix $J$ of the
  parametrization (\ref{reducedpara}) has full rank almost everywhere.
  The Jacobian matrix has the form
    \[
  J \,\,\, = \,\,\,
  \bordermatrix{
    & \sigma & \lambda\cr
    \psi_{ii} & I_p & B\cr
    \psi_{ij} & 0 & A
  }.
  \] 
  The entries in the unit matrix $I_p$ on the upper left are
  \begin{equation}
    \label{eq:dpsi}
  \frac{\partial \psi_{ij}}{\partial\sigma_{tt}} \,\, = \,\,
  \begin{cases}
    1 &\text{if}\; t=i=j,\\
    0 &\text{else}.
  \end{cases}
  \end{equation}
  The  matrix $J$ has full rank if and only if the
  $\binom{p}{2}\times \big(pm-\binom{m}{2}\big)$-matrix $A$ has full
  rank. The entries of the latter matrix $A$ are  
  \begin{equation}
    \label{eq:dlambda}
  \frac{\partial \psi_{ij}}{\partial\lambda_{st}} \,\, = \,\, 
  \begin{cases}
    \lambda_{tt} &\text{if}\; i<j\;\text{and}\; (i,j)=(t,s),\\
    \lambda_{jt} &\text{if}\; i=s=t<j\; \text{or}\; t<i=s<j, \\
    \lambda_{it} &\text{if}\; t<i<j=s,\\
    0 &\text{else}.
  \end{cases}
  \end{equation}
  If we set $\psi_{i<}=(\psi_{i,i+1},\dots,\psi_{i,p})^t \in\RRR^{p-i}$ and
  $\lambda_{>j}=(\lambda_{j+1,j},\dots,\lambda_{p,j})^t \in\RRR^{p-j}$,
  then the matrix $A$ can be written in the following form 
  \begin{equation*}
    \bordermatrix{ 
    & \lambda_{11} & \lambda_{22} & \lambda_{33}  
    & \!\!\dots\!\! &
    \lambda_{mm} & &  
    \lambda_{>1} & \lambda_{>2} & \lambda_{>3} & \!\!\dots\!\! &
    \lambda_{>m}\cr 
    \psi_{1<} & \lambda_{>1}&&&&&\vrule&A_{11}\cr 
    \psi_{2<} & &\lambda_{>2}&&&& \vrule&A_{21}&A_{22}\cr
    \psi_{3<} &&&\lambda_{>3}&&& \vrule&A_{31}&A_{32}&A_{33}\cr
    \vdots & &&&\ddots&& \vrule&\vdots&\vdots&\vdots&\ddots\cr    
    \psi_{m<} &&&&&\lambda_{>m}&
    \vspace{-0.2cm}
    \vrule&A_{m1}&A_{m2}&A_{m3}&\ddots&A_{mm}\cr
\cline{2-12}
    \psi_{m+1,<} &&&&&&\vrule&A_{m+1,1}& A_{m+1,2}&A_{m+1,3}&\dots&A_{m+1,m}\cr
    \vdots & &&&&& \vrule&\vdots&\vdots&\vdots&\cdots&\vdots\cr    
    \psi_{p-1,<} &&&&&& \vrule&A_{p-1,1}&A_{p-1,2}&A_{p-1,3}&\dots&A_{p-1,m}\cr
  }
  \end{equation*}
  \normalsize 
  where void entries are zero.  The submatrices
  $A_{ii}=\lambda_{ii}\, I_{p-i}$ are diagonal, and, for $i>j$,
  \[
  A_{ij} = 
  \bordermatrix{
    & \lambda_{j+1,j} &\dots & \lambda_{i-1,j} & \lambda_{i,j} &
    \lambda_{i+1,j} & \lambda_{i+2,j} & \dots&\lambda_{p,j}\cr
    \psi_{i,i+1} &&&&\lambda_{i+1,j}&\lambda_{ij}\cr
    \psi_{i,i+2} &&&&\lambda_{i+2,j}&&\lambda_{ij}\cr
    \vdots &&&&\vdots&&&\ddots\cr
    \psi_{i,p} &&&&\lambda_{p,j}&&&&\lambda_{ij}\cr
    }.
  \]
  Some of the submatrices in the partition of $A$ may not be present if $p$
  is too small.  For example, the entire lower half of $A$ is not present
  if $p\le m+1$.  If $p\le m+1$, then $A$ has a lower-triangular structure
  and is clearly of full rank if all $\lambda_{ii}$ are non-zero.  So we
  will assume that $p\ge m+2$, in which case the lower half of $A$
  comprises $\binom{p-m}{2}$ rows.
 
  We will now choose a particular matrix $\Lambda^0$ for which
  $A^0=A(\Lambda^0)$ is of full rank.  The existence of such $\Lambda^0$
  implies that the rank of $A$ is full for almost every choice of
  $\Lambda$.  The matrix $\Lambda^0$ has entries in $\{0,1\}$ with the
  non-zero entries chosen as follows.  For all $i\in [m]$, we set
  $\lambda^0_{ii}=1$.  As a consequence, the upper right block of $A$ is of
  full rank $\sum_{i=1}^m (p-i)=pm-\binom{m+1}{2}$.  The remaining non-zero
  entries of $\Lambda^0$ are determined as follows.  Let $J(p,m)$ be the
  minimum of $m$ and $\binom{p-m}{2}$. For $j\in [J(p,m)]$, let $i(j)$ be
  the integer in $\{m+1,\dots,p-1\}$ such that the $j$-th row of the lower
  half of $A$ is indexed by $\psi_{i(j),t}$ with $t\ge i(j)+1$.  For
  $j=1,\dots,J(p,m)$, we set exactly two components of the vector
  $\lambda_{>j}$ equal to one, namely those appearing in that row of
  $A_{i(j),j}$ that is part of the $j$-th row of the lower half of $A$.  As
  examples, consider $(p,m)=(7,4)$ and $(p,m)=(8,4)$, for which the above
  procedure selects the two (transposed) matrices
  \[
  (\Lambda_0)^t =
  \begin{pmatrix}
    1 & 0 & 0 & 0 & 1 & 1 & 0\\ 
    0 & 1 & 0 & 0 & 1 & 0 & 1\\
    0 & 0 & 1 & 0 & 0 & 1 & 1\\
    0 & 0 & 0 & 1 & 0 & 0 & 0
  \end{pmatrix}\quad
  \text{and}\quad
 (\Lambda_0)^t =
  \begin{pmatrix}
    1 & 0 & 0 & 0 & 1 & 1 & 0 & 0\\ 
    0 & 1 & 0 & 0 & 1 & 0 & 1 & 0\\
    0 & 0 & 1 & 0 & 1 & 0 & 0 & 1\\
    0 & 0 & 0 & 1 & 0 & 1 & 1 & 0
  \end{pmatrix}.
  \]
  
  Since the matrix $\Lambda^0$ has entries in $\{0,1\}$, the same holds for
  the matrix $A^0$.  A submatrix $A^0_{ij}$ of the upper right block of
  $A^0$, that is, $1\le j < i\le m$, has only one column with non-zero
  entries because $\lambda^0_{ij}=0$ if $i\le m$ and $i\not= j$.  This
  non-zero column is indexed by $\lambda_{ij}$ and can thus be eliminated
  by subtracting the row of $A^0$ indexed by $\psi_{ji}$.  This way the
  upper right block of $A^0$ is transformed into a unit matrix of size
  $\sum_{i=1}^m (p-i)=pm-\binom{m+1}{2}$, while no fill-in occurs in the
  upper left block of $A^0$.  
  
  Next we eliminate the lower right block of $A^0$ by subtraction of rows
  from the upper half. This elimination creates fill-in in the lower left
  block of $A^0$.  This fill-in is zero except for $J(p,m)$ many entries
  that are all equal to $-2$.  These non-zero entries occur in the
  positions $(j,j)$, $j=1,\dots,J(p,m)$, within the lower left block of
  $A^0$.  It follows that the rank of $A^0$ is equal to
  \[
  pm-\binom{m+1}{2}+J(p,m)=
  \min\left\{ pm-\binom{m}{2}, \binom{p}{2}\right\},
  \]
  which is the minimum of the number of rows and columns of $A^0$. Hence,
  $A^0$ is of full rank as we had claimed.  This concludes the proof of the
  stated formula for $\dim(F_{p,m})$.  The codimension is $\binom{p+1}{2}$
  minus $\dim(F_{p,m})$, and inequality (\ref{ineqSolveForp}) is gotten by
  solving $\binom{p-m}{2} > m$ for $p$.
\end{proof}

\begin{example} \label{ex:pentad}
Let us consider the case of two factors $m = 2$. The model
$F_{p,2}$ has positive codimension if and only if $p \geq 5$.
For $p = 5$, Theorem \ref{thm:dim} says that
$F_{5,2}$ has codimension $1$, so it 
is a hypersurface in the space of
symmetric $5 \times 5$-matrices. The hypersurface is defined
by the polynomial
\begin{eqnarray*}
f \,\, = & \,\,\,
\psi_{12} \psi_{13} \psi_{24} \psi_{35} \psi_{45}
- \psi_{12} \psi_{13} \psi_{25} \psi_{34} \psi_{45}
- \psi_{12} \psi_{14} \psi_{23} \psi_{35} \psi_{45}
+ \psi_{12} \psi_{14} \psi_{25} \psi_{34} \psi_{35} \\ &
+ \psi_{12} \psi_{15} \psi_{23} \psi_{34} \psi_{45}
- \psi_{12} \psi_{15} \psi_{24} \psi_{34} \psi_{35}
+ \psi_{13} \psi_{14} \psi_{23} \psi_{25} \psi_{45}
- \psi_{13} \psi_{14} \psi_{24} \psi_{25} \psi_{35} \\ &
- \psi_{13} \psi_{15} \psi_{23} \psi_{24} \psi_{45}
+ \psi_{13} \psi_{15} \psi_{24} \psi_{25} \psi_{34}
- \psi_{14} \psi_{15} \psi_{23} \psi_{25} \psi_{34}
+ \psi_{14} \psi_{15} \psi_{23} \psi_{24} \psi_{35}.
\end{eqnarray*}
This is the {\em pentad constraint} which was
first derived by \cite{kelley:35}. If $\Psi$ is the covariance
matrix of a distribution in the model $\mathbf{F}_{5,2}$
then $f(\Psi) = 0$, and the pentad $f$ is the unique
irreducible polynomial (up to scalar multiplication)
with this property. In the next section we discuss
the use of such invariants as test statistics, and in
the subsequent sections we derive higher invariants
using methods of computational algebra.

For $p = 4$, Theorem \ref{thm:dim} says that
$F_{4,2}$ has codimension $0$, so it is 
full-dimensional in the space of 
symmetric $4 \times 4$-matrices. The theorem does not
state that every positive definite matrix
$\Psi$ is in the model $F_{4,2}$. All it states is that
the decomposition of Proposition  \ref{prop:factordef},
\begin{equation}
\label{estimation42}
 \Psi \quad = \quad
\begin{pmatrix}
\sigma_{11} & 0 & 0  & 0 \\
0 & \sigma_{22} & 0 & 0   \\
0 & 0 & \sigma_{33} & 0 \\
0  & 0 & 0 & \sigma_{44} \\
\end{pmatrix}
+
\begin{pmatrix}
\lambda_{11} &  \lambda_{12} \\
\lambda_{21} &  \lambda_{22} \\
\lambda_{31} &  \lambda_{32} \\
\lambda_{41} &  \lambda_{42} 
\end{pmatrix}
\cdot
\begin{pmatrix}
\lambda_{11} & \lambda_{21} & \lambda_{31} & \lambda_{41} \\
\lambda_{12} & \lambda_{22} & \lambda_{32} & \lambda_{42} 
\end{pmatrix},
\end{equation}
imposes {\em no equality constraints} on the covariance matrix $\Psi$.  But
it does impose constraints in the form of inequations $f(\Psi) \not= 0$ and
inequalities $f(\Psi) \geq 0$.  We will discuss this issue in Section 4,
after the algebraic set-up of ideals has been introduced. Note that the
statistical problem of {\em parameter identification} corresponds to the
algebraic problem of solving the equations (\ref{estimation42}) for the
unknowns $\sigma_{ii},\lambda_{ij}$ when the $\psi_{ij}$ are given.
\end{example}


\section{Invariants as test statistics}
\label{sec:tetradsstats}

Let $\Psi\in\RRR^{p\times p}$ be a covariance matrix, that is, a positive
definite symmetric $p \times p$-matrix, and let $f$ be a polynomial in the entries $\psi_{ij}$
of $\Psi$.  We write $f(\Psi)$ for the evaluation of $f$ using the
numerical values of a particular matrix $\Psi$.  The polynomial $f$ is
called an {\em invariant\/} of the factor analysis model $\mathbf{F}_{p,m}$
if $f(\Psi)=0$ for all matrices $\Psi$ in the parameter space $F_{p,m}$.
Classical examples of invariants are the tetrad and pentad. If $f$ is an
invariant of $\mathbf{F}_{p,m}$ and $\Psi$ is a covariance matrix such that
$f(\Psi)\not=0$ then we can deduce that $\Psi\not\in F_{p,m}$.  This
suggests that model invariants can be used as statistics in tests of model
fit. We propose the following approach for putting this on a sound basis.

Assume we observe a sample $X_1,X_2,\ldots,X_N$ of independent random vectors
in $\RRR^p$ that are identically distributed according to the multivariate
normal distribution $\ND_p(\mu,\Psi)$ with mean vector $\mu\in\RRR^p$ and
positive definite $p\times p$-covariance matrix $\Psi$.  Let $\bar
X=\frac{1}{N} \sum_{k=1}^N X_k$ be the sample mean vector and
consider the sample covariance matrix
\[
S \, = \,(s_{ij})
\,\,\, = \,\,\, \frac{1}{N-1}\sum_{k=1}^N (X_k-\bar X)(X_k-\bar X)^t.
\]
  Moreover, let $f$ be an invariant of a
hypothesized factor analysis model $\mathbf{F}_{p,m}$.  The {\em sample
  invariant\/} $f(S)$ provides a consistent estimator of the true invariant
evaluation $f(\Psi)$.  The variance of $f(S)$, which we denote by
$\Var_\Psi[f(S)]$, can be derived by computing appropriate moments of the
Wishart distribution according to which the matrix $(N-1)\cdot S$ is
distributed; compare \citet[Sect.~3.4]{mardia:79} and \citet{wishart:28b}.
The variance $\Var_\Psi[f(S)]$ is a polynomial function of the true
covariance matrix $\Psi$. Replacing $\Psi$ by the sample covariance
matrix $S$ in this polynomial yields the estimator $\Var_S[f(S)]$. 
Using this estimator, we can define the {\em
  standardized sample invariant\/}
\begin{equation}
  \label{SSI}
Z_f \,\,\, = \,\,\, \frac{f(S)}{\sqrt{\Var_S[f(S)]}}.
\end{equation}

\begin{proposition}
  \label{prop:zf}
  Let $f$ be an invariant of the model $\mathbf{F}_{p,m}$,
  and let $\Psi$ be a covariance matrix such that $f$
  vanishes at $\Psi$ but its gradient vector $\nabla f $
  does not vanish at $\Psi$. Then,
  as the sample size $N$ tends to infinity,
  the standardized sample invariant $Z_f$ converges in distribution to a
  standard normal distribution:
  \[
  Z_f \,\,\, \longrightarrow_d \,\,\, \ND(0,1).
  \]
  \end{proposition}
\begin{proof}
  The vectorization of $\sqrt{N}(S-\Psi)$ converges in distribution to a
  centered multivariate normal distribution. Hence, by the delta method
  \citep[p.279]{shorack:00}, 
  \[
  \sqrt{N}\, f(S) \,\,\, = \,\,\, \sqrt{N}\, [f(S)-f(\Psi)]
  \,\,\, \longrightarrow_d \,\,\, \ND(0,v_f)
  \]
  with asymptotic variance 
  \[v_f = \lim_{N\to \infty} N\Var_\Psi[f(S)] >0.
  \]  
  Since $\Var_\Psi[f(S)]/\Var_S[f(S)]$ converges in probability to one, it
  follows from Slutsky's theorem that
  \[
  Z_f = \frac{\sqrt{\Var_\Psi[f(S)]}}{\sqrt{ \Var_S[f(S)]}}\cdot
  \frac{\sqrt{N}\, f(S)}{\sqrt{N \Var_\Psi[f(S)]}} \longrightarrow_d
  \tfrac{1}{\sqrt{v_f}} \cdot \ND(0,v_f) = \ND(0,1). \qedhere
  \]
\end{proof}

\begin{remark}
  \label{rem:biascorrect}
  The sample invariant $f(S)$ is typically a biased estimator of $f(\Psi)$.
  However, if the expectation of $f(S)$ is of the form
  $\E_\Psi[f(S)]=h(N)\cdot f(\Psi)$, where $h(N)$ is a function of the sample
  size only, then one can consider the {\em bias-corrected sample
    invariant\/} $\tilde f(S) = f(S)/h(N)$.  An analog of Proposition
  \ref{prop:zf} holds when $f(S)$ is replaced by $\tilde f(S)$.
\end{remark}

\begin{example} \label{ex:tetradtest}
  We derive the standardized and bias-corrected sample invariants for
  the one-factor model $\mathbf{F}_{p,1}$. Let $i$, $j$, $k$, $\ell$
  be four distinct indices in $[p]$ and consider the tetrad
  \begin{equation}
    \label{eq:tetradgeneral}
    f \,\,\, = \,\,\, \psi_{ik}\psi_{j\ell}-\psi_{i\ell}\psi_{jk}.
  \end{equation}
  If $\Psi$ is a covariance matrix in $F_{p,1}$ then the tetrad
  vanishes, i.e., $f(\Psi)=0$.
  
  The sample tetrad $f(S)=s_{ik}s_{j\ell}-s_{i\ell}s_{jk}$ is a
  consistent but biased estimator of $f(\Psi)$.  However, the bias can
  be corrected as described in Remark \ref{rem:biascorrect} with the
  {\em bias-corrected sample tetrad\/} being equal to
  \[  
  \tilde f(S) \,\,\, = \,\,\,
  \frac{N-1}{N-2}\left(s_{ik}s_{j\ell}-s_{i\ell}s_{jk}\right). 
  \]
  For any covariance matrix $\Psi$, the variance of this unbiased estimator
  of $f(\Psi)$ is equal to
  \begin{equation}
    \begin{split}
      \label{eq:wishart}
      \Var_\Psi\big[ \tilde f(S)\big] &=
      \frac{N+1}{(N-1)(N-2)}\, \det(\Psi_{\{i,j\}\times\{i,j\}})\cdot
      \det(\Psi_{\{k,\ell\}\times\{k,\ell\}}) \\
      &\qquad- \frac{1}{N-2}
      \det(\Psi_{\{i,j,k,\ell\}\times\{i,j,k,\ell\}})
      + \frac{3}{N-2}
      \det(\Psi_{\{i,j\}\times\{k,\ell\}})^2.
          \end{split}
  \end{equation}
  This expression was first computed by \cite{wishart:28}.  If
  $\Psi\in F_{p,1}$, then
  $$ 
  \det(\Psi_{\{i,j\}\times\{k,\ell\}}) \,\,\, = \,\,\,
  \psi_{ik}\psi_{j\ell}-\psi_{i\ell}\psi_{jk}
  \,\,\, = \,\,\, 0.
  $$
  Thus the last term in (\ref{eq:wishart}) vanishes and we can use
  the estimate
  \begin{multline*}
  \Var_S\big[ \tilde f(S)\big] =\\
  \frac{N+1}{(N-1)(N-2)}\, \det(S_{\{i,j\}\times\{i,j\}})\cdot
  \det(S_{\{k,\ell\}\times\{k,\ell\}}) - \frac{1}{N-2}
  \det(S_{\{i,j,k,\ell\}\times\{i,j,k,\ell\}}).
  \end{multline*}
  Following the recipe in (\ref{SSI}), we 
 introduce the standardized bias-corrected sample tetrad 
  \[
  Z_{\tilde f} \,\,\, = \,\,\, \frac{\tilde f(S)}{\sqrt{\Var_S\big[
        \tilde f(S)\big]}}.
  \] 
  This is an explicit expression which can be evaluated for any sample
  covariance matrix $S$ arising from data $X_i$.  If at least one of
  the four entries $\psi_{ik}$, $\psi_{j\ell}$, $\psi_{i\ell}$,
  $\psi_{jk}$ is non-zero, then the gradient of the tetrad is non-zero
  at $\Psi$.  Proposition \ref{prop:zf} says that $Z_{\tilde f}$ has
  an asymptotic standard normal distribution $\ND(0,1)$ when the
  sample size $N$ tends to infinity.
\end{example}  
  
Suppose now that $f$ is an arbitrary polynomial invariant of the
factor analysis model $\mathbf{F}_{p,m}$, and we wish to test the null
hypothesis $H_f:f(\Psi)=0$.  In light of Proposition \ref{prop:zf}, we
can do this by computing the corresponding standardized sample
invariant $Z_f$ and by comparing it to the appropriate quantile of the
standard normal distribution $\ND(0,1)$.  More precisely, for chosen
significance level $\alpha\in (0,1)$, we can find an interval
$[-c_\alpha,c_\alpha]$ which, assuming $H_f$ is true, contains the
standardized sample invariant $Z_f$ with (asymptotic) probability
$1-\alpha$.  If we observe a value $z_f$ of $Z_f$ that falls outside
this interval, $z_f\not\in [-c_\alpha,c_\alpha]$, then this
constitutes evidence against $H_f$ and, in particular, evidence
against the hypothesized factor analysis model of which $f$ is an
invariant.
  
If the hypothesized factor analysis model is a hypersurface then
consideration of a single invariant $f$ is sufficient.  This happens,
for instance, in the case $p=5,m=2$ discussed in Example
\ref{ex:pentad}. Here $f$ is the pentad and we only need to test
$H_f$.
  
In general, however, the model structure will not be captured in a single
polynomial invariant.  Then we might want to employ a set $I$ of invariants
of the considered model to test model fit.  For instance, $I$ could be a
set of ideal generators as in Section 4.  A simple approach to working with
several invariants is to employ Bonferroni's inequality, which suggests the
consideration of the interval $[-c_{\alpha/|I|},c_{\alpha/|I|}]$.  This
interval simultaneously contains all standardized sample invariants $Z_f$,
$f\in I$, with probability at least $1-\alpha$.  Therefore, if one or more
observed values $z_f$ fall outside the interval
$[-c_{\alpha/|I|},c_{\alpha/|I|}]$, then we have found statistical evidence
against the hypothesized factor analysis model.  More powerful approaches
than this simple Bonferroni method can be obtained by combining the
invariants in a quadratic form; see e.g.~\cite{hipp:03} for work
employing tetrads.  Alternatively, \cite{spirtes:00} employ tests of
vanishing tetrads to define scores for model selection in Gaussian
graphical models with hidden variables.
  
Tetrads appear to be the only invariants of factor analysis models that
have seen routine use in data analysis.  However, the approach we have
outlined above is feasible also for other invariants such as the pentad and
the higher invariants we determine subsequently.  The only difficulty
involved is the estimation of the variance-covariance structure of the
sample invariants.  We expect that recent work on moments of the Wishart
distribution \citep{massam:05,richards:01} can be applied fruitfully to
overcome this difficulty.  When moments of invariants cannot be determined
exactly, asymptotic approximations can be derived from the asymptotic
covariance matrix for the sample covariance matrix; see \cite{roverato:98}
for a discussion of properties of the Isserlis matrix which determines this
asymptotic covariance matrix.


\section{Algebraic setup}
\label{sec:algsetup}

We are interested in polynomial relations among the entries of a factor
analysis covariance matrix $\Psi\in F_{p,m}$.  The mathematical framework
for studying such polynomial relations is that of commutative algebra and
algebraic geometry \citep[see e.g.][]{CLS:97}.  Many algorithms from these
fields are implemented in software for symbolic computation, and they
provide powerful computational tools for the study of model invariants.
The application of these tools in statistics is the focus of
algebraic statistics \citep{pistone:01,ascb}. While algebraic statistics
has so far been predominantly occupied with the study of models for
discrete random variables, the present study is one of the first in
this emerging field which concerns continuous random variables. The set-up
to be introduced is fairly general and can be used to study arbitrary
Gaussian graphical models, not just the factor analysis model.

We fix the ring of polynomials with real coefficients in the
$\binom{p+1}{2}$ indeterminates $\psi_{ij}$:
\[
\begin{split}
  \RRR[ \psi_{ij},\; i\le j   ] \quad = \quad
  \RRR[\psi_{11},\psi_{12},\ldots,\psi_{1p},
  \psi_{22},\ldots,\psi_{2p}, \psi_{33},\ldots,\psi_{pp}].
\end{split}
\]
For any subset $F$ of the symmetric matrices in $\RRR^{p\times p}$, let
$I(F)$ be the set of all polynomials $f\in \RRR[\psi_{ij},\; i\le j]$ such
that $f(\Psi)=0$ for all $\Psi\in F$.  Clearly, $I(F)$ is an ideal in
$\RRR[\psi_{ij},\; i\le j]$; that is, the sum of two polynomials in $I(F)$
is again in $I(F)$, and the product of any polynomial in $\RRR[\psi_{ij},\;
i\le j]$ with a polynomial in $I(F)$ is in $I(F)$.  According to Hilbert's
basis theorem, every ideal is generated by a finite list of polynomials. We
tacitly assume this finite representation for all the ideals which appear
in the following discussion.

The object of our interest is the {\em ideal of invariants\/} of the model
$\mathbf{F}_{p,m}$, which is the ideal $I_{p,m}=I(F_{p,m})$.  Since
membership in $F_{p,m}$ depends only on the off-diagonal entries of the
matrix $\Psi$, we can regard $I_{p,m}$ as an ideal in the subring
$\,\RRR[\psi_{ij},\; i <j]\,$ of $\,\RRR[\psi_{ij},\; i\le j]$.  If $I$ is
any ideal in the bigger polynomial ring $\RRR[\psi_{ij},\; i\le j]$ then
the intersection $\, I \,\cap \, \RRR[\psi_{ij},\; i <j]\,$ is an ideal in
the smaller polynomial ring $\RRR[\psi_{ij},\; i <j]$.  Passing to this
intersection is the process of {\em elimination} of the variables
$\psi_{11}, \ldots,\psi_{pp}$.  The following result shows how the ideal
$I_{p,m}$ can be computed using elimination.

\begin{theorem}\label{thm:elim}
  Let $M_{p,m}\subseteq \RRR[\psi_{ij},\; i\le j]$ be the ideal that
  is generated by all $(m+1)\times (m+1)$-minors of a symmetric matrix
  $\Psi\in\RRR^{p\times p}$.  Then the ideal of invariants equals
\begin{equation}
\label{thm4id}
  I_{p,m} \quad = \quad M_{p,m} \,\cap \, \RRR[\psi_{ij},\; i < j].
\end{equation}
\end{theorem}

\begin{proof}
  The proof makes use of standard arguments from algebraic geometry, and
  all varieties $V( \,\,\cdot \,\,)$ are understood over the field $\CCC$
  of complex numbers.  Recall that the {\em variety} $V(M_{p,m})$ of the
  ideal $M_{p,m}$ is the set of common zeroes of the polynomials in
  $M_{p,m}$. This set coincides with the set of all symmetric $p \times
  p$-matrices of rank at most $m$.  Let $F'_{p,m}$ denote the set of all $p
  \times p$-matrices of the form $\Psi = \Sigma + \Gamma$ where $\Sigma$ is
  a diagonal matrix and $\Gamma \in V(M_{p,m})$. Thus $F'_{p,m}$ is the
  superset of our parameter space $F_{p,m}$ gotten by dropping the positive
  definiteness requirement.  Since the cone of positive definite matrices
  is open (and hence Zariski dense) in the space of all symmetric matrices,
  we conclude that $F_{p,m}$ and $F'_{p,m}$ have the same Zariski closure
  in $\CCC^{p \times p}$.  The projection of this Zariski closure onto the
  space of off-diagonal entries is Zariski closed, and it coincides with
  the variety $V(I_{p,m}) $ of the desired ideal $I_{p,m}$.  On the other
  hand, every matrix in the projection of $F'_{p,m}$ is the image of a
  matrix $\Gamma$ in $V(M_{p,m})$. We conclude that $V(I_{p,m})$ equals
  the Zariski closure of the projection of $V(M_{p,m})$ onto the
  off-diagonal coordinates.
  
  
  By the Elimination Theorem \citep[see][\S 3.2]{CLS:97} we have
  \begin{equation}
    \label{SameVarieties}
    V \bigl(M_{p,m} \,\cap \, \RRR[\psi_{ij},\; i < j] \bigr) \quad = \quad
    V(I_{p,m}).  
  \end{equation}
  Now, it is known that $M_{p,m}$ is a prime ideal and that the minors form
  a Gr\"obner basis for $M_{p,m}$; compare \cite{conca:94}, \citet[][Corollary
  4.11, Example 4.12]{ss:05}.  The primality of $M_{p,m}$ implies that
  $\,M_{p,m} \cap \RRR[\psi_{ij},\:i < j ]\,$ is prime as well.  Since, by
  definition, $I_{p,m}$ is radical, we apply Hilbert's Nullstellensatz to
  (\ref{SameVarieties}) to conclude that $\,I_{p,m} \, = \, M_{p,m} \,\cap
  \, \RRR[\psi_{ij},\; i < j]$.
\end{proof}

Theorem \ref{thm:elim} allows for the derivation of a finite
generating set of the ideal $I_{p,m}$ using the method of Gr\"obner
bases. This will be explained in Section 5. In the remainder of this
section, we discuss some consequences and geometric aspects of Theorem
\ref{thm:elim}, starting with some polynomials that obviously belong
to the ideal $I_{p,m}$.  Up to sign, there are
\[
\frac{1}{2}\binom{p}{2(m+1)}\binom{2(m+1)}{m+1} 
\]
 {\em off-diagonal $(m+1)\times (m+1)$-minors} of the matrix $\Psi$,
that is, subdeterminants that do not involve any diagonal entries of $\Psi$.
Such minors of size $m+1$ are trivially in $I_{p,m}$.

\begin{corollary}
  \label{cor:minor}
  Let $p\ge 2(m+1)$ and choose two disjoint sets $R,C\subset [p]$ of
  cardinality $|R|=|C|=m+1$.  Then the off-diagonal minor $\det(\Psi_{R\times
    C})$ is in $I_{p,m}$.
\end{corollary}

\begin{example}
  Let $m = 1$.  Then the $2 \times 2$-off-diagonal minors of $\Psi$ belong
  to $I_{p,1}$.  For example, if $R = \{1,2\}$ and $C = \{3,4\}$ then this minor is the tetrad
  $$
  \det(\Psi_{R\times C}) \,\,\, = \,\,\, \det \begin{pmatrix}
      \psi_{13} & \psi_{14} \\
      \psi_{23} & \psi_{24}
\end{pmatrix}  \,\,\, = \,\,\,
\psi_{13}\psi_{24} - \psi_{14} \psi_{23} .  $$
\end{example}

If $p<2(m+1)$ there are no off-diagonal minors.  However, if $p\ge 2m+1$,
it is still easy to determine some non-zero polynomials in $I_{p,m}$ by
considering two minors that contain exactly one common diagonal entry
$\psi_{ii}$ and eliminating this diagonal entry.

\begin{corollary}
  \label{prop:fi}
  Let $m\ge 2$ and $p\ge 2m+1$, and choose $i\in [p]$.  Let $(R,C)$ and
  $(\bar R,\bar C)$ be two pairs of $m$-element subsets of $[p] \backslash
  \{i\}$ that are disjoint: $R\cap C=\emptyset=\bar R\cap \bar C$.
  Finally, let $\Psi^0$ be the symmetric matrix whose off-diagonal entries
  are the unknowns $\psi_{ij}$ and whose diagonal entries are equal to 0.
  Then the following polynomial is in $I_{p,m}$:
    \begin{equation*}
    f_{i,R,C,\bar R,\bar C} \quad = \quad 
    \det(\Psi_{R\times C})\cdot 
    \det(\Psi^0_{(i,\bar R)\times (i,\bar C)})
    -  
    \det(\Psi_{\bar R\times \bar C})\cdot 
    \det(\Psi^0_{(i,R)\times  (i,C)}) .
  \end{equation*}
  The notation $(i,R)\times (i,C)$ indicates that the $i$-th row and column
  are arranged first.
\end{corollary}

\begin{proof} 
  The polynomial $ f_{i,R,C,\bar R,\bar C}$ lies in $\RRR[\psi_{ij},\; i <
  j]$.  The two identities 
  \[
  \begin{split}
    \det(\Psi_{(i,R)\times (i,C)}) \quad &= \quad
    \det(\Psi^0_{(i,R)\times (i,C)}) + \,\psi_{ii} \cdot \det(\Psi_{R\times
      C}), \\
    \det(\Psi_{(i,\bar R)\times (i,\bar C)}) \quad &= \quad
    \det(\Psi^0_{(i,\bar R)\times (i,\bar C)}) + \,\psi_{ii} \cdot
    \det(\Psi_{\bar R\times \bar C})
  \end{split}
  \] 
  imply
  \begin{equation*}
    f_{i,R,C,\bar R,\bar C} \quad = \quad 
    \det(\Psi_{R\times C})\cdot 
    \det(\Psi_{(i,\bar R)\times (i,\bar C)})
    -  
    \det(\Psi_{\bar R\times \bar C})\cdot 
    \det(\Psi_{(i,R)\times  (i,C)}).
  \end{equation*}
  This is a polynomial linear combination of $(m+1)\times
  (m+1)$-minors of $\Psi$, so it lies in $M_{p,m}$.  We conclude that
  $f_{i,R,C,\bar R,\bar C}$ is in the right hand side of (\ref{thm4id})
  and hence in $I_{p,m}$.
\end{proof}

 The \emph{linear eliminant} $\,f_{i,R,C,\bar R,\bar C} \,$ is a homogeneous
 polynomial of degree $2m+1$.  For $m = 2$ and $m = 3$, the linear
 eliminants recover the tetrads and the pentads as follows:

\begin{remark}
  \label{rem:tets}
  If $m=1$ and $p\ge 4$ then we can choose pairs $(R,C)$ and $(\bar R,\bar
  C)$ satisfying the assumptions of Corollary \ref{prop:fi}.  The result is
  a polynomial combination of two tetrads:
  \[
  f_{i,R,C,\bar R,\bar C} \,= \, -\psi_{rc}\psi_{i\bar r}\psi_{i\bar c}
  + \psi_{\bar r\bar c}\psi_{ir}\psi_{ic} \, = \, -\psi_{i\bar
    c}\cdot\det(\Psi_{\{r,\bar r\}\times \{c,i\}}) +
  \psi_{ir}\cdot\det(\Psi_{\{\bar c,c\}\times\{\bar r,i\}}).
  \]
\end{remark}

\begin{example}
  Let $m =2$ and $p = 5$.  Then the polynomial $f_{i,R,C,\bar R,\bar C}$
  has degree five, and it does not depend on the choices of $i$, $R$, $C$,
  $\bar{R}$ and $\bar{C}$.  Up to sign, it coincides with the pentad $f$
  which was displayed in Example \ref{ex:pentad}.  Note that the twelve
  monomials in the pentad $f$ correspond to the twelve labeled cycles on
  the set of nodes $\{1,2,3,4,5\}$.  The ideal $I_{5,2}$ is the principal
  ideal generated by the pentad; in symbols, $I_{5,2} = \langle f \rangle
  $.
\end{example}

The following proposition shows that linear eliminants
are non-redundant invariants.

\begin{proposition}\label{prop:tight}
  Let $m\ge 2$ and $p\ge 2m+1$.  If $R\cup C=\bar R\cup\bar C$, then linear
  eliminant $f_{i,R,C,\bar R,\bar C}$ is not in the ideal generated by the
  off-diagonal $(m \! + \! 1)\times (m\! + \! 1)$-minors of $\Psi$.
  \end{proposition}

\begin{proof}
  Without loss of generality assume that $i=1$ and $R\cup
  C=\{2,3,\dots,2m+1\}$.  Since $\codim(F_{2m+1,m})>0$, by Theorem
  \ref{thm:dim}, we can choose a symmetric matrix $\Psi \in\RRR^{p\times
    p}$ such that (i) $f_{1,R,C,\bar R,\bar C}(\Psi)\not=0$ and (ii) all
  off-diagonal entries in row $2m+2$ to $p$ are zero.  Then all
  off-diagonal $(m+1)\times (m+1)$-minors of the chosen matrix $\Psi$ are
  zero. This shows that $f_{1,R,C,\bar R,\bar C}$ cannot be a polynomial
  combination of the off-diagonal minors.
\end{proof}

We believe that the following converse to Proposition \ref{prop:tight} holds.   
As we will see, Conjecture \ref{conj:tight} is part of a general series of
finiteness conjectures about the ideals $I_{p,m}$.

\begin{conjecture}\label{conj:tight}
  Let $m \ge 2$ and $p \ge 2m+2$.  If $R \cup C \neq \bar{R} \cup \bar{C}$
  then the linear eliminant $f_{i,R,C,\bar R,\bar C}$ is in the ideal
  generated by the off-diagonal $(m \! + \! 1) \times (m \! + \! 1)$-minors
  of $\Psi$.
\end{conjecture}
We call the linear eliminants where $R \cup C = \bar{R} \cup \bar{C}$
the {\em $(2m +1)$-ads}.  So when $m = 2$, we recover the pentads and
when $m =3$ we obtain {\em septads}.

We close this section with a discussion
of the geometric role played by the ideal $I_{p,m}$
in the context of factor analysis.
Recall that the variety $V(I_{p,m})$
is the set of all common zeros of the
polynomials in the ideal $I_{p,m}$.
Consider the following four statements:

\begin{itemize}
\item[(a)]
A polynomial vanishes on the parameter space $F_{p,m}$
if and only if it lies in $ I_{p,m}$.
\item[(b)]
The factor analysis model $\mathbf{F}_{p,m}$
is represented by the variety $V(I_{p,m})$.
\item[(c)]
The parameter space $F_{p,m}$
coincides with the variety $V(I_{p,m})$.
\item[(d)]
The closure of the parameter space $F_{p,m}$
coincides with the variety  $V(I_{p,m})$.
\end{itemize}

Then statement (a) is true because this is how the
ideal $I_{p,m}$ was defined. Statement (b) is vague, but it
expresses the philosophy of this paper, so we simply
declare it to be true. On the other hand, statement (c) is false 
with respect to every meaningful interpretation of what
the statement may mean. In algebraic geometry,
$V(I_{p,m})$ denotes the set of zeros of $I_{p,m}$ over the field $\CCC$
of complex numbers, and this is what was
meant in the proof of Theorem \ref{thm:elim}.
Considering the zeros of $I_{p,m}$
among positive-definite matrices, positive semi-definite matrices,
or just real symmetric matrices, we get the  inclusions
\[
V_{\rm pd}(I_{p,m}) \subset V_{\rm psd} (I_{p,m}) 
\subset V_\RRR(I_{p,m}) \subset V(I_{p,m}). 
\]
So, meaningful interpretations of (c) may be that $F_{p,m}$ equals
$V_{\rm pd}(I_{p,m})$, and that $V(I_{p,m})$ equals the set $F^\CCC_{p,m}$ of
complex $p \times p$-matrices $\Psi = \Sigma + \Lambda \Lambda^t $ where
$\Sigma$ is diagonal and $\Lambda \in \CCC^{p \times m}$.  Both of these
statements are false as the following example shows.

\begin{example}
  Let $p = 3$ and $m = 1$. Then $F^\CCC_{3,1}$ consists of all $3
  \times 3$-matrices of the form
  \begin{equation}
    \label{estimation31}
    \begin{pmatrix}
      \psi_{11} & \psi_{12} & \psi_{13} \\
      \psi_{12} & \psi_{22} & \psi_{23} \\
      \psi_{13} & \psi_{23} & \psi_{33}
    \end{pmatrix}
    \quad = \quad
    \begin{pmatrix}
      \sigma_{11} & 0 & 0   \\
      0 & \sigma_{22} & 0    \\
      0 & 0 & \sigma_{33}  \\
    \end{pmatrix}
    +
    \begin{pmatrix}
      \lambda_{11}  \\
      \lambda_{21}  \\
      \lambda_{31}
    \end{pmatrix}
    \cdot
    \begin{pmatrix}
      \lambda_{11} & \lambda_{21} & \lambda_{31}
    \end{pmatrix}.
  \end{equation}
  The three off-diagonal identities imply
  \begin{equation}
    \label{eq:offdiag} 
    \lambda_{11}^2 \psi_{23} - \psi_{12} \psi_{13} \,\,\, = \,\,\,
    \lambda_{21}^2 \psi_{13} - \psi_{12} \psi_{23} \,\,\, = \,\,\,
    \lambda_{31}^2 \psi_{12} - \psi_{13} \psi_{23}   \quad = \quad 0 .
  \end{equation}
  This shows that a matrix in $\Psi$ in $F_{3,1}^\CCC$ cannot have
  precisely one zero off-diagonal entry.  But $V(I_{3,1})$ consists of
  all symmetric $3 \times 3$-matrices since ${\rm codim}(F_{3,1}) = 0$
  and $I_{3,1} = \{0\}$.  Hence $F_{3,1}^\CCC$ is a proper subset of
  $V(I_{3,1})$, and, likewise, $F_{3,1}$ is a proper subset of $V_{\rm
    pd}(I_{3,1})$.
\end{example}

Let us now come to statement (d). This statement is true over the field
$\CCC$ of complex numbers. Every matrix in $V(I_{p,m})$ is the limit of
matrices in $F_{p,m}^\CCC$.  This follows from a non-trivial algebraic
geometry result to the effect that, for the image of any polynomial map
over $\CCC$, the usual closure coincides with the Zariski closure 
\citep[see][Proposition 7, p.~490]{CLS:97}.
 On the other hand, statement (d) is false over the
real numbers.  Namely, in our example, $V_{\rm pd}(I_{3,1})$ is the set of
all positive-definite matrices. If $\Psi$ is a positive-definite matrix
with $\psi_{12} > 0$, $\psi_{13} > 0$ and $\psi_{23} < 0$, then
(\ref{eq:offdiag}) forces $\lambda_{11}$ to be the square root of a
negative number, so $\Psi$ cannot be in the closure of $F_{3,1}$.  A
similar (but more complicated) analysis can be performed for the case
$p=4,m=2$ starting from the equations given in (\ref{estimation42}).

In summary, in this paper we do not determine all the constraints satisfied
by the parameter space $F_{p,m}$ of the factor analysis model.  What we do
determine is the set $I_{p,m}$ of all polynomial equation constraints.
These characterize the closure of $F_{p,m}$ if we allow complex numbers.
The polynomials in $I_{p,m}$ are the model invariants, and, as argued in
Section \ref{sec:tetradsstats}, they can be used to derive novel test statistics for Gaussian
graphical models.


\section{Gr\"obner basis computations}
\label{sec:ideals}

We now focus on computing finite generating sets for the ideals
$I_{p,m}$.  Following Theorem \ref{thm:elim}, this can be done by
equating to zero all $(m\!+\!1) \times (m\!+\!1)$-minors of an unknown
symmetric $p \times p$-matrix $\Psi$, and then eliminating the
off-diagonal unknowns $\psi_{ii}$ from these equations. In computer
algebra, there are two main methods for eliminating unknowns from a
system of equations: Gr\"obner bases and resultants.  In this section
we present the Gr\"obner basis approach, while resultants will be
featured in the next section.  We shall assume familiarity with
``Gr\"obner basics'' at the level of \cite{CLS:97}.

The complete answer to our problem is currently only known for the
one-factor model $(m=1)$. Namely, as we shall see in Theorem
\ref{thm:tetradgb}, the tetrads provide a reduced Gr\"obner basis for
the ideal $I_{p,1}$. For two or more factors $(m \geq 2)$, we did
numerous computations with the computer algebra systems {\tt
  Macaulay2} and {\tt Singular}. The results of these computations are
presented in this section (see Tables \ref{tab:codimdegree} and
\ref{tab:betti} below).  We shall return to these results in Section
\ref{sec:conjectures}, where we offer some conjectures about the
ideals $ I_{p,m}$.

Let $m=1$. For any four indices $i<j<k<\ell$ in $[p]$, we have the
tetrads in (\ref{eq:tetrad}). Since the first tetrad is the difference
of the third and the second tetrad, it suffices to pick out the last
two tetrads in (\ref{eq:tetrad}).  Let
\[
\mathcal{T}_p \quad = \quad \{
\underline{\psi_{ij}\psi_{k\ell}}-\psi_{ik}\psi_{j\ell},\:
\underline{\psi_{i\ell}\psi_{jk}}-\psi_{ik}\psi_{j\ell}\mid 1\le
i<j<k<\ell\le p\}
\]
be the set of $2\binom{p}{4}$ tetrads obtained in this way.  As
described by \citet{deloera:1995}, the underlined terms are the
leading terms with respect to a certain monomial order $\succ$ on
$\RRR[\psi_{ij},\; i< j]$.  (They call this monomial order the {\em
  thrackle order}.)

\begin{theorem}
  \label{thm:tetradgb}
  If $p\le 3$ the ideal $I_{p,1}$ is the zero ideal.  If $p\ge 4$, the
  set $\mathcal{T}_p$ is the reduced Gr\"obner basis of the ideal
  $I_{p,1}$ with respect to the monomial order $\succ$.
\end{theorem}
\begin{proof}
  The claim follows from Theorem 2.1 in \cite{deloera:1995}.
\end{proof}

If we observe $p=5$ variables, then the set $\mathcal{T}_5$ contains
ten tetrads.  In \citet[p.~76]{harman:76} it is stated that one can
find five of these tetrads such that {\em ``any other conditions must be
linearly dependent on the [five tetrads].'}'  Similarly,
\citet[p.~277]{hipp:03} state that {\em ``to detect the full set of
redundant vanishing tetrads when there are more than four variables
requires careful algebraic derivation''}  and {\em ``in the case of the
five-indicator model there will be five nonredundant vanishing
tetrads.'' } \citet[p.~418]{harman:76} outlines a justification of his
claim, which, however, is valid only if all $\psi_{ij}$
are non-zero.  In a strict algebraic sense, Harman's claim is
incorrect because none of the
ten tetrads in $\mathcal{T}_5$ is a polynomial
linear combination of the other nine tetrads. 

\vskip .1cm

Moving on to the general case $m \geq 2$, we 
now demonstrate how to compute a minimal generating set
of the ideal $I_{p,m}$ by means of two software
packages for algebraic geometry.

\begin{example}[$p=7$, $m=2$ in {\tt Macaulay 2}]
\label{ex:m2}
The first software we used is the program {\tt Macaulay 2} due
to \cite{m2:98}.
To compute a minimal generating set of the ideal 
$I_{7,2}$ using {\tt Macaulay 2},
we use the following sequence of six commands:

\begin{small}
\begin{verbatim}
R = QQ[p11,p22,p33,p44,p55,p66,p77,p12,p13,p14,p15,p16,p17,p23,p24,p25,p26,p27,
       p34,p35,p36,p37,p45,p46,p47,p56,p57,p67, MonomialOrder=>Eliminate 7];
Psi = matrix{{p11,p12,p13,p14,p15,p16,p17},
             {p12,p22,p23,p24,p25,p26,p27},
             {p13,p23,p33,p34,p35,p36,p37},
             {p14,p24,p34,p44,p45,p46,p47},
             {p15,p25,p35,p45,p55,p56,p57},
             {p16,p26,p36,p46,p56,p66,p67},
             {p17,p27,p37,p47,p57,p67,p77}};
M72 = minors(3,Psi);
I72 = ideal selectInSubring(1,gens gb M72);
mingens I72
codim I72, degree I72
\end{verbatim}
\end{small}

\noindent
The command {\tt I72 = ideal selectInSubring(1,gens gb M72)} performs
the actual elimination step of deriving $\,I_{7,2} \, $ from $\, M_{7,2}$. The command
{\tt mingens I72} outputs $56$ polynomials which minimally
generate the ideal $I_{7,2}$. This list includes $35$ polynomials of degree
three and $21$ polynomials of degree five.  The latter are  $21$ pentads like

\begin{small}
\begin{verbatim}
 p36p37p45p47p56-p35p37p46p47p56-p36p37p45p46p57+p35p36p46p47p57
+p34p37p46p56p57-p34p36p47p56p57+p35p37p45p46p67-p35p36p45p47p67
-p34p37p45p56p67+p34p35p47p56p67+p34p36p45p57p67-p34p35p46p57p67.
\end{verbatim}
\end{small}

\noindent
Of the $35$ polynomials of degree three, $21$ are off-diagonal
minors like 

\begin{small}
\begin{verbatim}
p26p35p47-p25p36p47-p26p34p57+p24p36p57+p25p34p67-p24p35p67.
\end{verbatim}
\end{small}

\noindent
The remaining $14$ polynomials are sums of off-diagonal minors. In fact,
we can replace these by $14$ off-diagonal minors such that
the resulting $35$ polynomials are minimal generators of
$I_{7,2}$.  This validates the entry for $p=7,m=2$  in Table \ref{tab:betti}.
Finally, the last command line informs us that the variety
$V(I_{7,2})$ has codimension $8$ and degree $259$.
\end{example}

The notion of ``degree'' requires an explanation.
Next to the codimension, this is the most important invariant
of an algebraic variety. Suppose that $V$ is a variety
of codimension $c$ in $\CCC^r$. Then the {\em degree} of $V$
is the number of points in $\,V \,\cap L \,$
where $L$ is a general affine subspace of dimension $c$ in $\CCC^r$.
The case $c =1 $ is familiar: if $V$ is a hypersurface, defined by
the vanishing of one polynomial $f$, then the degree of $V$
equals the degree of $f$, and this 
is the number of intersection points of $V$ with a general line.

Table \ref{tab:codimdegree} summarizes what we know about the 
codimension and the degree of the factor analysis model $V(I_{p,m})$
for $m \leq 5$ and $p \leq 9$. In this section we discuss
the $m=2$ and $m=3$ columns, and in Section 6
we discuss $(m,p)=(4,8)$ and $(m,p)=(5,9)$.

\begin{table}[htbp]
  \centering
  \begin{tabular}{r|rr|rr|rr|rr|rr}
     & \multicolumn{2}{c|}{$m=1$} & \multicolumn{2}{c|}{$m=2$} &
    \multicolumn{2}{c|}{$m=3$} & \multicolumn{2}{c|}{$m=4$} &
    \multicolumn{2}{c}{$m=5$}\\  
    $p$ & \text{codim} & \text{deg} & \text{codim} & \text{deg}
    & \text{codim} & \text{deg} & \text{codim} & \text{deg} & \text{codim}
    & \text{deg} \\  
    \hline
    3 & 0  &   1 &  0 &   1  &  0 &    1 &  0 &  1 & 0 &  1\\
    4 & 2  &   4 &  0 &   1  &  0 &    1 &  0 &  1 & 0 &  1\\
    5 & 5  &  11 &  1 &   5  &  0 &    1 &  0 &  1 & 0 &  1\\
    6 & 9  &  26 &  4 &  45  &  0 &    1 &  0 &  1 & 0 &  1\\
    7 & 14 &  57 &  8 & 259  &  3 &   91 &  0 &  1 & 0 &  1\\
    8 & 20 & 120 & 13 & 1232 &  7 & 1368 &  2 & 98 & 0 &  1\\
    9 & 27 & 247 & 19 & 5319 & 12 &14232 &  6 &  ? & 1 & 54 
  \end{tabular}
  \caption{Codimensions and degrees for the factor analysis model.}
  \label{tab:codimdegree}
\end{table}

For $m = 2$ our computations suggest that the ideal $I_{p,2}$ is always
generated by minors and pentads, but so far we have been unable to prove
this for general $p$. See Section \ref{sec:conjectures} for a discussion of
this conjecture. For $m \geq 3$ we found that the computer algebra system
{\tt Singular} performs better than {\tt Macaulay 2}.  Here is a
non-trivial computation.

\begin{example}[$p=8$, $m=3$ in {\tt Singular}] \label{ex:p8m3}
  To compute the ideal $I_{8,3}$ using {\tt Singular} \citep{GPS01}, the
  following sequence of eight commands can be used:

\begin{small}
\begin{verbatim}
ring R = 0,(p11,p22,p33,p44,p55,p66,p77,p88,
            p12,p23,p34,p45,p56,p67,p78,p18,
            p13,p24,p35,p46,p57,p68,p17,p28,
            p14,p25,p36,p47,p58,p16,p27,p38, p15,p26,p37,p48),dp;
matrix Psi[8][8] = p11,p12,p13,p14,p15,p16,p17,p18,
                   p12,p22,p23,p24,p25,p26,p27,p28,
                   p13,p23,p33,p34,p35,p36,p37,p38,
                   p14,p24,p34,p44,p45,p46,p47,p48,
                   p15,p25,p35,p45,p55,p56,p57,p58,
                   p16,p26,p36,p46,p56,p66,p67,p68,
                   p17,p27,p37,p47,p57,p67,p77,p78,
                   p18,p28,p38,p48,p58,p68,p78,p88;
ideal M83 = minor(Psi,4);
ideal I83 = eliminate(M83, p11*p22*p33*p44*p55*p66*p77*p88);
nvars(basering) - dim(I83); mult(I83); // codimension and degree
ideal I83min = mstd(I83)[2];  // minimal generators
betti(I83min);
I83min;
\end{verbatim}
\end{small}

\noindent
In the first command, which declares the polynomial ring, we list the
variables in an order different from the order chosen in the {\tt Macaulay
  2} code described in Example \ref{ex:m2}.  Together with the option {\tt
  dp}, this variable ordering determines a monomial order which we found to
be advantageous.  With this particular order, a modern workstation requires
less than 10 minutes to execute the remaining seven commands. Other
monomial orders we considered led to significantly slower computations.

The function {\tt eliminate} carries out the elimination step of deriving
$I_{8,3}$ from $M_{8,3}$. The codimension and degree of the variety
$V(I_{8,3})$ are equal to 7 and 1368, respectively, as computed in the line
following the elimination.  The function {\tt mstd} permits to compute the
list of polynomials {\tt I83min}, which minimally generate the ideal
$I_{8,3}$.  According to the command {\tt betti(I83min)}, this
generating set consists 
of $14$ polynomials of degree four, $260$ polynomials of
degree seven, and $168$ polynomials of degree eight. They are:

\begin{description}
\item[\rm\em degree 4:] Twelve polynomials are off-diagonal $4\times
  4$-minors.  The other two polynomials are sums of off-diagonal minors and
  can be replaced by off-diagonal minors.

\item[\rm\em degree 7:] The 260 polynomials of degree seven come in two flavors.
  \begin{enumerate}
  \item[(a)] Of the 260 polynomials, $120=\binom{8}{7}\cdot 15$ are
    sums of linear eliminants of the type described in Proposition
    \ref{prop:tight}.  These linear eliminants are equal to
    $f_{i,R,C,\bar R,\bar C}$ with $R\cup C=\bar R\cup \bar C$ in
    $[p]\setminus\{i\}$ and we call them {\em septads}.  A septad has
    168 terms, and an example of a septad appearing in the minimal
    generator is

    \begin{small}
\begin{verbatim}
p23p13p46p68p17p28p47+p67p18p13p24p28p36p47-p13p24p68p17p28p36p47-...
...+p12p17p36p26p37p48^2-p23p17p16p26p37p48^2+p13p16p27p26p37p48^2 .
\end{verbatim}
    \end{small}

    \noindent
    While 102 of the 120 polynomials in our output are 
    septads, the other twelve can be replaced by septads without affecting the
    minimal generator property.
  \item[(b)] The remaining $260-120=140$ polynomials are of ideal-theoretic
    nature.  They are not polynomial linear combinations of the off-diagonal
    minors and septads. However, the squares of the
     $140$ polynomials are polynomial combinations of off-diagonal minors
    and septads.  Hence these invariants vanish at a covariance matrix 
    $\Psi$ whenever the off-diagonal minors and septads do.
      \end{enumerate}
\item[\rm\em degree 8:] These $168$ minimal generators are also of ideal-theoretic
  nature.  Their squares are polynomial linear combinations of off-diagonal minors
  and septads.
\end{description}
\end{example}

Table \ref{tab:betti} summarizes our knowledge about the composition of a
minimal generating set of the factor analysis ideal $I_{p,m}$ for $m\le 3$ and $p\le 9$.  
This table suggests various conjectures
about the ideals $I_{p,m}$ for general $p$, and we will discuss these in
Section \ref{sec:conjectures}.

\begin{table}[htbp]
  \centering
  \begin{tabular}{c|r|rr|rrrr}
     & \multicolumn{1}{c|}{$m=1$} & \multicolumn{2}{c|}{$m=2$} & 
    \multicolumn{4}{c}{$m=3$}\\
    $p$&deg 2& deg 3& deg 5 & deg 4& deg 7& deg 7& deg 8\\
    \hline
    4&   2& --- & ---& --- & --- & --- & ---\\
    5&  10&  0 &   1 & --- & --- & --- & ---\\
    6&  30&  5 &   6 & --- & --- & --- & ---\\
    7&  70& 35 &  21 &  0 & 15   &   0 &  20\\
    8& 140& 140 & 56 & 14 & 120 &  140 & 168\\
    9& 252& 420 & 126 & 126 &  540 &  1386 & 756\\
    \hline
    & {\em tetrad} & {\em minor} & {\em pentad} & {\em minor} & {\em septad} &
    {\em ideal-} & {\em ideal-}\\
    &&&&&& {\em theor.} & {\em theor.} 
  \end{tabular}
  \caption{The degrees of the minimal generators of the ideals $I_{p,m}$.}
  \label{tab:betti}
\end{table} 


\section{Multilinear resultants}
\label{sec:resultants}

Resultants are a technique for simultaneously eliminating
$m$ unknowns from a system of $m+1$ polynomial equations.
While Gr\"obner bases can be used to perform this task,
we must remember that Gr\"obner bases are a very general method.
They often compute too much and are hence too inefficient.
The applicability of resultants is more limited, but
wherever they do apply, resultants tend to outperform
Gr\"obner bases. 
For general introductions to resultants
see \cite{gkz:94} and \cite{sturmfels:97,sturmfels:02}.

In our algebraic study of factor analysis, we found the Gr\"obner
basis computations of Section 5 to be infeasible for $m \geq 4$.
Instead we did some computations using the {\em multilinear
  resultant}. In this section we explain this technique, and how it
was used to derive the degrees $54$ and $98$ for $(m,p)=(5,9)$ and
$(m,p)=(4,8)$ in Table \ref{tab:codimdegree}.

 Consider a set of $n+1$ multilinear polynomials
$f_0,\ldots,f_n$  in $n$ unknowns $x_1,\ldots,x_n$:
\begin{equation}
\label{m+1INm}
 f_j \quad  = \quad \!\!\!\!\!
\sum_{i_1,i_2,\ldots,i_n \in \{0,1\}} \!\!\!
  a^{j}_{i_1 i_2 \cdots i_n} x_1^{i_1} x_2^{i_2} \cdots x_n^{i_n}
\qquad \qquad (j=0,1,\ldots,n).
\end{equation}
Here the coefficients $a^j_{i_1 i_2 \cdots i_n}$ are regarded as unknowns.
The total number of these coefficients is $2^n \cdot (n+1)$, and they
generate a polynomial ring which we denote by
$$ \RRR [{\mathbf a}] \; := \;
\RRR \bigl[\, a^{j}_{i_1 \cdots i_n} \,: \,\,
i_1,\ldots,i_n \in \{0,1\}, j \in \{0,\ldots,n\} \bigr]. $$
We write $\RRR[\mathbf{a}, \mathbf{x}]$ for the
polynomial ring generated by the coefficients
$a^j_{i_1 i_2 \cdots i_n}$ and the unknowns
$x_1,\ldots,x_n$, and we consider the ideal
$\,\langle f_0,f_1,\ldots,f_n \rangle\,$  in $\RRR[ \mathbf{a}, \mathbf{x}]$
which is generated by the multilinear polynomials (\ref{m+1INm}).
We have the following result from algebra:

\begin{theorem}
The elimination ideal $\,\langle f_0,f_1,\ldots,f_n \rangle\,\cap \,
 \RRR [\mathbf{a}] \,$ is generated by an irreducible
 polynomial ${\mathcal R}({\bf a})$ which is homogeneous of
 degree $\,n \, ! \,$ in the coefficients of each $f_j$.
\end{theorem}

\begin{proof}
  This follows from the results in \citet[Sect.~8.2.A]{gkz:94}, applied to
  the special case when the toric variety $X_A$ is the product of $n$
  projective lines in its Segre embedding.  The corresponding polytope $Q$ is
  the $n$-dimensional standard cube, which has normalized volume equal to
  $\,n! \, $.  The polynomial ${\mathcal R}({\bf a})$ is the Chow form of
  $X_A$.
\end{proof}

We call ${\mathcal R}({\bf a})$ the {\em $n$-th multilinear resultant}.
Here are the first three cases:

\begin{example}
\label{resultantn=1}
If $n = 1$ then $f_0  = a^0_0 + a^0_1 x_1 $ and
$f_1 = a^1_0 + a^1_1 x_1 $. Their resultant equals
\begin{equation}
\label{n=1RES}
{\mathcal R}({\bf a}) \,\,\, = \,\,\,
a^0_0 a^1_1 - a^0_1 a^1_0 \,\,\, = \,\,\,
a^1_1 \cdot f_0 - a^0_1 \cdot f_1 \,\,\, = \,\,\, 
\left| \begin{array}{cc}
a^0_0 & a^0_1 \\
a^1_0 & a^1_1 
\end{array} \right|. 
\end{equation}
\end{example}

\begin{example}
\label{resultantn=2}
If $n=2$ then we are considering a system of three bilinear equations
\begin{eqnarray*}
f_0 \,\,\, = & a^0_{00} + a^0_{10} x_1 + a^0_{01} x_2 + a^0_{11} x_1 x_2 ,\\
f_1 \,\,\, = & a^1_{00} + a^1_{10} x_1 + a^1_{01} x_2 + a^1_{11} x_1 x_2 ,\\
f_2 \,\,\, = & a^2_{00} + a^2_{10} x_1 + a^2_{01} x_2 + a^2_{11} x_1 x_2 .
\end{eqnarray*}
Their resultant has the following determinantal representation:
\begin{equation}
\label{n=2RES}
{\mathcal R}({\bf a}) \, = \,
\left| \begin{array}{ccc}
\! a^0_{00} & a^0_{10} & a^0_{01}  \! \\
\! a^1_{00} & a^1_{10} & a^1_{01} \! \\
\! a^2_{00} & a^2_{10} & a^2_{01}  \!
\end{array} \right|
 \cdot
\left| \begin{array}{ccc}
 \! a^0_{10} & a^0_{01} & a^0_{11} \! \\
 \! a^1_{10} & a^1_{01} & a^1_{11} \! \\
 \! a^2_{10} & a^2_{01} & a^2_{11} \!
\end{array} \right| \,-\,
\left| \begin{array}{ccc}
\! a^0_{00} &  a^0_{01} & a^0_{11} \! \\
\! a^1_{00} &  a^1_{01} & a^1_{11} \! \\
\! a^2_{00} &  a^2_{01} & a^2_{11} \!
\end{array} \right|
\cdot
\left| \begin{array}{ccc}
\! a^0_{00} & a^0_{10} &  a^0_{11} \! \\
\! a^1_{00} & a^1_{10} &  a^1_{11} \! \\
\! a^2_{00} & a^2_{10} &  a^2_{11} \!
\end{array} \right| .
\end{equation}
This polynomial has degree six but it is quadratic
in the coefficients of each $f_j$.
\end{example}

\begin{example}
\label{resultantn=3}
If $n=3$ then the coefficients of $f_0,f_1,f_2,f_3$ 
form an $4 \times 8$-matrix
\begin{equation}
\label{Cmatrix}
A \quad = \quad \begin{bmatrix}
a^0_{000} & a^0_{001} &  a^0_{010} & a^0_{011} &
a^0_{100} & a^0_{101} & a^0_{110} & a^0_{111} \\
a^1_{000} & a^1_{001} &  a^1_{010} & a^1_{011} &
a^1_{100} & a^1_{101} & a^1_{110} & a^1_{111} \\
a^2_{000} & a^2_{001} &  a^2_{010} & a^2_{011} &
a^2_{100} & a^2_{101} & a^2_{110} & a^2_{111} \\
a^3_{000} & a^3_{001} &  a^3_{010} & a^3_{011} &
a^3_{100} & a^3_{101} & a^3_{110} & a^3_{111} \\
\end{bmatrix}.
\end{equation}
Let $[ijkl]$ denote the determinant of the $4 \times 4$-submatrix of $A$
with columns $i,j,k,l$.  Then the multilinear resultant $ {\mathcal R}({\bf
  a}) $ is the determinant of the following $6 \times 6$-matrix:

\begin{small}
\begin{equation}
\label{n=3RES}
\begin{bmatrix}
 [0124] & [0234] & [0146]-[0245] &
    [0346]-[0247] & [0456] & [0467] \\
  &  &  & & & \\
 [0125] & [1234] &
   [0147]+[0156]&
   -[1247]+[0356]&
   [1456] & [1467] \\
 +[0134]  & +[0235] & -[0345]-[1245]  & 
   -[0257]+[1346] & +[0457] & +[0567] \\
  &  &  & & & \\
   [0135] & [1235] & [0157]-[1345] &
    -[1257]+[1356] & [1457] & [1567] \\
    &  &  & & & \\
  [0126]  & [0236] & -[1246]+[0256] &
   [2346]-[0267] &  [2456] & [2467] \\
   &  &  & & & \\
  [0136] & [1236] &
  -[1247]-[1346] &
  -[0367]-[1267] &
  [3456] & [2567] \\
+[0127]  & +[0237] & +[0257]+[0356]  & +[2356]+[2347] & +[2457]&
+[3467] \\
   &  &  & & & \\
  [0137] & [1237] & -[1347]+[0357] &
  -[1367]+[2357] & [3457] &  [3567] \\
\end{bmatrix}
\end{equation}
\end{small}

\noindent
The derivation of this $6 \times 6$-matrix is explained in
\citet[Proposition 4.10]{sturmfels:02}. The matrix (\ref{n=3RES}) differs
from the $6 \times 6$-matrix displayed by \citet{sturmfels:02} because the
latter had some typographical errors. These typos have now been corrected
in (\ref{n=3RES}).
\end{example}

These formulas can be applied to algebraic factor analysis as follows.
Recall from Theorem \ref{thm:elim} that the ideal $I_{p,m}$ is computed by
eliminating the diagonal unknowns $\psi_{ii}$ from the $(m+1) \times
(m+1)$-minors of an indeterminate covariance matrix $\Psi = (\psi_{ij})$.
Multilinear resultants are relevant for this task because the minors are
multilinear polynomials in the $\psi_{ii}$. For instance, the derivation of
the pentad constraint can be interpreted as an evaluation of the first
multilinear resultant (\ref{n=1RES}). Likewise, the second multilinear
resultant (\ref{n=2RES}) can be used to produce non-trivial invariants in
$I_{p,m}$ when $p \geq 2, \, m \geq 8$.

The general method for computing invariants of the factor analysis model
using resultants is the following.  We
choose any subset $D=\{d_1,\dots,d_n\}$ of cardinality $n$ in $[p]$ and any
subsets $R_0, \ldots,R_n, C_0,\ldots,C_n\,$ each of cardinality $m+1-n$ in
$[p]\backslash D$ such that $\,R_i \,\cap\,C_i \,= \,\emptyset \,$ for all
$i$.  This choice of $D,R_\bullet,C_\bullet$ gives rise to an invariant as
follows.

\begin{theorem} \label{resultantdegree}
  Let $x_j = \psi_{d_j d_j}$ for $j=1,\ldots,n$ and
  $\,f_k = {\rm det}(\Psi_{D R_k \times D C_k})\,$ for $k = 0,\ldots,n$.
  The evaluation of the $n$-th multilinear resultant at $f_0,f_1,\ldots,f_n$ is an
  element of $\, \RRR[ \psi_{ij},\; i\le j ]$.  This polynomial lies in the
  ideal $I_{p,m}$, and it is either zero or it is homogeneous of degree
\begin{equation}
\label{resdegreeformula}
 ( m \, - \,\frac{n}{2}\,+ \,1) \cdot
(n+1)\, !
\end{equation}
\end{theorem}

\begin{proof}
The proof  uses the
methods of \cite{gkz:94} and is omitted here.
\end{proof}

It is instructive to examine this theorem for small values of $n$.
If $n=1$ (as in Example \ref{resultantn=1})
then the degree (\ref{resdegreeformula}) equals $2m+1$
and we recover the linear eliminant
$\,f_{i,R,C,\bar R,\bar C} \,$
of Proposition \ref{prop:fi}.
Here
$D = \{i\}$,
$R_0 = R$,
$R_1 = \bar R$,
$C_0 = C$ and $
C_1 = \bar C$.
If $n=2$ (as in Example \ref{resultantn=2}) then
the invariant constructed in Theorem \ref{resultantdegree}
is homogeneous of degree $6m$.
If $n=3$ (as in Example \ref{resultantn=3}) then
the invariant constructed in Theorem \ref{resultantdegree}
is homogeneous of degree $24m-12$. In particular, the
degree (\ref{resdegreeformula}) is $108$ for $m=5$.

\begin{example}[$p=9$, $m=5$ using resultants]
\label{ex:res95}
This is the case appearing in the lower 
right corner of Table \ref{tab:codimdegree}.
The variety $V(I_{9,5}) $ is a hypersurface of degree $54$
in the space of symmetric $9 \times 9$-matrices.
The irreducible polynomial defining
this hypersurface is the greatest common divisor (gcd) of 
all multilinear resultants constructed for  $n=3$
as in Theorem \ref{resultantdegree}.
Each resultant has degree $108$, and their gcd was found
to have degree $54$. In fact, we observed that it
 suffices to take the gcd  of only two such resultants.

For a concrete instance take $D = \{1,2,3\}$,
$R_0 = \{4,5,6\}$, $C_0 = \{7,8,9\}$,
$R_1 = \{4,5,7\}$, $C_1 = \{6,8,9\}$,
$R_2 = \{4,6,8 \}$, $C_2 = \{5,7,9 \}$,
$R_3 = \{4,7,9\}$, $C_3 = \{5,6,8\}$,
and let $f$ be the degree $108$ polynomial
constructed by Theorem \ref{resultantdegree}.
It is infeasible to express $f$ as a
sum of monomials. (Note that the number of
monomials of degree $108$ in the $36$ unknowns
$\psi_{ij}$ exceeds $10^{33}$). However, 
the $6 \times 6$-determinant
(\ref{n=3RES}) offers an efficient 
representation of the invariant $f$.
Namely, if the $\psi_{ij}$ are replaced 
by linear forms in one or two parameters, then
this  $6 \times 6$-determinant can 
be evaluated rapidly. From this
we see that $f$ is non-zero and that it 
factors as the product of four irreducible
polynomials, having degrees $18$, $18$, $18$ and $54$.
The last factor is the generator of $I_{9,5}$.

Let $g$ denote the degree $54$ invariant which generates the ideal
$I_{9,5}$. We have represented $g$ as the gcd of several polynomials of
degree $108$, but this representation specifies $g$ only up to a
multiplicative constant.  Up to this constant, we can evaluate $g$
numerically at a covariance matrix $\Psi$ by choosing a random matrix
$\Psi_0$, introducing an unknown, say $t$, computing the gcd of several
resultants of the form $f(\Psi+t\cdot \Psi_0)\in\RRR[t]$, and evaluating the resulting polynomials at $t = 0$.  Repeating the
computation at a second covariance matrix $\Psi'$ with the same $\Psi_0$ is an efficient scheme
for evaluating the ratio $g(\Psi)/g(\Psi')$.
\end{example}

\begin{example}[$p=8$, $m=4$ using {\tt Macaulay 2} and resultants]
\label{ex:res84}
The variety $V(I_{8,4}) $ has codimension two
in the space of symmetric $8 \times 8$-matrices.
In order to compute its degree,
we intersect the projective variety of
$I_{8,4} $ with a general two-dimensional plane.
We do this computation over a finite field,
using the following {\tt Macaulay 2}  commands:

\begin{small}
\begin{verbatim}
S = ZZ/101[r,s,t];
R = ZZ/101[p11,p22,p33,p44,p55,p66,p77,p88,r,s,t,MonomialOrder=> Eliminate 8];
f = map(R,S,{r,s,t});
    p12 = f(random(1,S));   p13 = f(random(1,S));   p14 = f(random(1,S)); 
    p15 = f(random(1,S));   p16 = f(random(1,S));   p17 = f(random(1,S)); 
    p18 = f(random(1,S));   p23 = f(random(1,S));   p24 = f(random(1,S)); 
    p25 = f(random(1,S));   p26 = f(random(1,S));   p27 = f(random(1,S));
    p28 = f(random(1,S));   p34 = f(random(1,S));   p35 = f(random(1,S)); 
    p36 = f(random(1,S));   p37 = f(random(1,S));   p38 = f(random(1,S)); 
    p45 = f(random(1,S));   p46 = f(random(1,S));   p47 = f(random(1,S)); 
    p48 = f(random(1,S));   p56 = f(random(1,S));   p57 = f(random(1,S));
    p58 = f(random(1,S));   p67 = f(random(1,S));   p68 = f(random(1,S)); 
                                                    p78 = f(random(1,S)); 
Psi = matrix {{p11,p12,p13,p14,p15,p16,p17,p18},
              {p12,p22,p23,p24,p25,p26,p27,p28},
              {p13,p23,p33,p34,p35,p36,p37,p38},
              {p14,p24,p34,p44,p45,p46,p47,p48},
              {p15,p25,p35,p45,p55,p56,p57,p58},
              {p16,p26,p36,p46,p56,p66,p67,p68},
              {p17,p27,p37,p47,p57,p67,p77,p78},
              {p18,p28,p38,p48,p58,p68,p78,p88}}
G = gens gb minors(5,Psi);   J = ideal selectInSubring(1,G);
codim J, degree J
\end{verbatim}
\end{small}

The output of repeated runs verifies that
the degree of $V(I_{8,4})$ equals $98$. However,
this computation gives no information about the
minimal generators of $I_{8,4}$, and at present
we do not even know the smallest degree of a non-zero
polynomial in $I_{8,4}$. 

On the other hand, we obtain a large number of non-trivial invariants
of degree $24$ by applying Theorem \ref{resultantdegree} with $n=2 $.
Namely, for any choice of indices $D, R_\bullet, C_\bullet$ we set
$$
f_i(x_1,x_2) \quad = \quad
{\rm det} \begin{pmatrix}
x_1 & \psi_{d_1,d_2 } &
 \psi_{d_1, c_{i0} } & 
 \psi_{d_1, c_{i1}} & 
 \psi_{d_1 ,c_{i2}} \\
 \psi_{d_2, d_1} & x_2 &
 \psi_{d_2, c_{i0}} & \psi_{d_2,c_{i1} } 
& \psi_{d_1 ,c_{i2}} \\
 \psi_{r_{i0},d_1 } &
 \psi_{r_{i0},d_2 } &
 \psi_{r_{i0},c_{i0} } &
 \psi_{r_{i0},c_{i1} } &
 \psi_{r_{i0},c_{i2} } \\
 \psi_{r_{i1},d_1 } &
 \psi_{r_{i1},d_2 } &
 \psi_{r_{i1},c_{i0} } &
 \psi_{r_{i1},c_{i1} } &
 \psi_{r_{i1},c_{i2} } \\
 \psi_{r_{i2},d_1 } &
 \psi_{r_{i2},d_2 } &
 \psi_{r_{i2},c_{i0} } &
 \psi_{r_{i2},c_{i1} } &
 \psi_{r_{i2},c_{i2} } 
\end{pmatrix}
\qquad \hbox{for} \,\,\, i=0,1,2.
$$
Here, $R_i=\{r_{i0},r_{i1},r_{i2}\}$ and $C_i=\{c_{i0},c_{i1},c_{i2}\}$.
The invariant of degree $24$ is obtained evaluating the
formula (\ref{n=2RES}), in which we abbreviate the coefficients of $f_i$ by
$\,a^i_{jk} = a^i_{jk}(\psi)$.  Just as in Example
\ref{ex:res95}, it is impossible to write this invariant as a sum of
monomials, but it is very easy to evaluate it numerically using the
determinantal representation (\ref{n=2RES}).
\end{example}


\section{Conjectures about generators of the ideals $I_{p,m}$}
\label{sec:conjectures}

In Sections 4-6 we discussed the problem of computing a finite
generating set for the ideal $I_{p,m}$ of invariants of the
factor analysis model $\mathbf{F}_{p,m}$.
Our computational results suggest
some natural conjectures and problems about the structure 
of this generating set, and we believe that these will be
of independent interest to commutative algebraists.

 The pattern we found for small $m$, and that we hope is true for larger $m$, is that as $p$ gets large the generators of $I_{p,m}$ depend on only a certain fixed number of random variables.  This type of \emph{finiteness} property has frequent occurrences in algebraic statistics. See  \cite{allman:04} and \cite{santos:03} for two instances.

The prototypical conjecture of this type is the following.

\begin{conjecture}\label{conj:2factor}
The ideal  of the two-factor model, $I_{p,2}$,
 is minimally generated by $5 \binom{p}{6}$ off-diagonal
$3 \times 3$-minors and $\binom{p}{5}$ pentads.
\end{conjecture}

Conjecture \ref{conj:2factor} is supported by the numerical evidence
compiled in Table \ref{tab:betti}.  Note that all of the minors and pentads
described involve at most six random variables.  However, unlike in the
case of the 1-factor model, there does not seem to be any term order that
makes the collection of off-diagonal minors and pentads a Gr\"obner basis
for $I_{p,2}$.

The most natural term order we discovered in our computations has already
been introduced in Example \ref{ex:p8m3}.  In general, this is the
lexicographic term order with $\psi_{ij} \succ \psi_{kl}$ if the circular
distance between $i$ and $j$ is smaller than the circular distance between
$k$ and $l$.  (To compute the circular distance between $i$ and $j$, place
the numbers $1, \ldots, p$ equispaced around a circle and measure the
distance by taking the shortest path around the circle between $i$ and
$j$.)  If these circular distances are the same, we declare $\psi_{ij}
\succ \psi_{kl}$ if $i < k$.  We call this term order the \emph{circular
  lexicographic term order}.  The fact that each of the ideals $I_{p,m}$ is
an $m$-th secant ideal, together with the machinery developed in
\cite{ss:05}, and our computations led us to the following conjecture.

\begin{conjecture}\label{conj:delight}
  The circular lexicographic term order is $2$-delightful for $I_{p,1}$.
  More specifically, the reduced Gr\"obner for $I_{p,2}$ consists of
  certain explicitly constructed polynomials of odd degree less than $p$
  with squarefree initial terms.
\end{conjecture}

The notion of a delightful term order was developed in \cite{ss:05} and
the technical details are beyond the scope of this short section.  However,
the basic idea is that, if Conjecture \ref{conj:delight} is true,
information about the initial ideal and reduced Gr\"obner basis of the
two-factor ideal $I_{p,2}$ can be deduced from the reduced Gr\"obner basis
of the one-factor ideal $I_{p,1}$ using graph theory.  We refer the
interested reader to \cite{ss:05} for information about delightful term
orders.

Moving to three factors, the situation seems even more complicated.
As we have already seen in Table \ref{tab:betti}, the minors and
septads are not enough to generate the entire ideal $I_{p,3}$.  The
minors and septads do, however, determine the parameter space
$F_{p,3}$ set-theoretically in all the examples we were able to
compute.  This leads us to suspect:

\begin{conjecture}
  Let $J_{p,3}$ be the ideal generated by all the $4 \times 4$ off-diagonal
  minors and all the septad linear eliminants. Then the radical of
  $J_{p,3}$ is the prime ideal $I_{p,3}$.
\end{conjecture}

We do not have any concrete conjectures about the generators or Gr\"obner
bases of $I_{p,3}$.  Note that each of the minors and septads involve at
most eight random variables.  Is it possible that this type of finiteness
behavior continues for larger $m$?  For a subset $A \subset [p]$ denote by
$I_{A,m}$ the ideal $I_{|A|,m}$ with indices labeled by the elements of
$A$.

\begin{question}\label{ques:finite}
For each  integer $m  \geq 1$ does there exist another integer $s(m)$ such that 
$$
I_{p,m} \,\,\,\, = \sum_{A \subset [p] , |A| = s(m)} \!\!\!\!\!\!
I_{A,m} \qquad \hbox{for all} \,\, p > s(m) \,\, ? 
$$
Does there exist a different number $t(m)$ where the equality holds up
to radical?
\end{question}

Our computations aside, there is some theoretical evidence to suggest that
the set-theoretic finiteness result might hold with $t(m) = 2m + 2$.
Namely, we can show that the set-theoretic finiteness result does hold in
the complement of a certain hypersurface.

\begin{proposition}\label{prop:2m+2}
  Let $\Psi$ be a symmetric $p\times p$-matrix with $p \geq 2m+2$ and
  suppose that no $m \times m$-minor of $\Psi$ is zero.  Then $\,\Psi \in
  F_{p,m}\,$ if and only if $\,\Psi_{A,A} \in F_{2m+2,m}$ for all subsets
  $A \subset [p]$ of cardinality $|A| = 2m+2$.
\end{proposition}

\begin{proof}
  The ``only if'' direction is trivial.  To prove the ``if'' direction we
  first need to refer to two simple results about factor analysis models.
  First of all, if $\Psi \in F_{p,m}$ and has no $m \times m$ minor equal
  to zero, then the decomposition $\Psi = \Sigma + \Gamma$ with $\Sigma$
  diagonal and ${\rm rank}(\Gamma) = m$ is unique.  Furthermore, if
  $\Gamma$ is a symmetric rank $m$ matrix and $\Lambda$ and $K$ are $p
  \times m$ matrices with $\Gamma = \Lambda\Lambda^t = KK^t$ then there
  exists an orthogonal matrix $Q$ such that $\Lambda = KQ$.  See
  \cite{anderson:56} for both of these results.

  Now we prove the ``if'' direction by induction on $p$.  The induction
  base is $p = 2m+2$ in which case the statement is vacuous.  Suppose that
  $\Psi$ satisfies $\Psi_{A,A} \in F_{2m+2,m}$ for all subsets $A \subset
  [p]$ with $|A| = 2m+2$.  Denote by $\Psi^+$ the submatrix $\Psi_{[p-1],
    [p-1]}$, by $\Psi_-$ the submatrix $\Psi_{[p] \backslash \{1\}, [p]
    \backslash \{1\}}$ and by $\Psi^+_-$ the submatrix $\Psi_{[p]
    \backslash \{1,p\}, [p] \backslash \{1,p\}}$.  In other words, $\Psi^+$
  is the upper left $(p-1) \times (p-1)$ submatrix, $\Psi_-$ is the lower
  right $(p-1) \times (p-1)$ submatrix, and $\Psi^+_-$ is the $(p-2) \times
  (p-2)$ submatrix where $\Psi^+$ and $\Psi_-$ overlap.

  By the induction hypothesis, $\Psi^+$ and $\Psi_-$ belong to $F_{p-1,m}$,
  and so there exist unique $\Sigma^+$, $\Sigma_-$ diagonal and $\Gamma^+$,
  $\Gamma_-$ of rank $m$ such that $\Psi^+ = \Sigma^+ + \Gamma^+$ and
  $\Psi_- = \Sigma_- + \Gamma_-$.  Furthermore, these provide a unique
  representation for the overlap $\Psi^+_- = \Sigma^+_- + \Gamma^+_-$,
  where $\Sigma^+_-$ and $\Gamma^+_-$ are the common overlapping portions
  of $\Sigma^+$ and $\Sigma_-$ and, respectively, $\Gamma^+$ and
  $\Gamma_-$.  Let $\Gamma^+ = \Lambda\Lambda^t$ and $\Gamma_- = KK^t$ be
  rank $m$ factorizations of $\Lambda^+$ and $\Lambda_-$.  Because of the
  uniqueness in the overlap $\Lambda^+_-$, and the second result taken from
  \cite{anderson:56} mentioned above, we can assume that the last $p-2$
  rows of $\Lambda$ coincide with the first $p-2$ rows of $K$.  Now form
  the $p \times p$ matrices
  $$\bar \Sigma = \left( \begin{array}{cc}
      \Sigma^+_{11} & 0 \\
      0 & \Sigma_- \end{array} \right) \quad \mbox{and} \quad \bar \Gamma =
  \bar{\Lambda} \bar{\Lambda}^t \quad \mbox{where} \quad \bar{\Lambda} =
  \left( \begin{array}{c}
      \Lambda_1 \\
      K
    \end{array} \right).
  $$
  Set $\bar \Psi = \bar \Sigma + \bar \Gamma$.  We claim that $\bar \Psi
  = \Psi$ and hence $\Psi \in F_{p,m}$, which completes the proof.  Note
  that trivially $\bar \psi_{ij} = \psi_{ij}$ except possibly for the pair
  $(i,j) = (1,p)$.  So we must show that $\bar \psi_{1p} = \psi_{1p}$.  Let
  $B$ and $C$ be disjoint subsets of $[p] \backslash \{1,p\}$ of
  cardinality $|B| = |C| = m$.  Then the following three $(m+1) \times
  (m+1)$ minors are equal to zero:
  $$
  \left| \begin{array}{cc}
      \bar \Psi_{1 \times C} & \bar \psi_{1p} \\
      \bar \Psi_{B \times C} & \bar \Psi_{B \times p} \end{array} \right| =
  \left| \begin{array}{cc}
      \Psi_{1 \times C} & \bar \psi_{1p} \\
      \Psi_{B \times C} & \Psi_{B \times p} \end{array} \right| = \left|
    \begin{array}{cc}
      \Psi_{1 \times C} & \psi_{1p} \\
      \Psi_{B \times C} & \Psi_{B \times p} \end{array} \right| = 0.
  $$
  The first and last minors are zero because $\bar \Psi_{A \times A}$
  and $\Psi_{A \times A}$ belong to $F_{2m+2,p}$ for any $A$ with $|A| =
  2m+2$ by assumption.  The middle minor is zero because it is entry-wise
  equal to the first minor.  Since by assumption ${\rm det} ( \Psi_{B
    \times C}) \neq 0$ we deduce that $\psi_{1p} = \bar \psi_{1p}$.
\end{proof}

The proof technique we present does not allow the extra condition in
Proposition \ref{prop:2m+2} to be dropped and thus we are a long way from
answering Question \ref{ques:finite}.


\section{Computer algebra for Gaussian models: the next steps}
\label{sec:future}

The research presented in this paper merely scratches the surface of
possible applications of computer algebra techniques for studying the
factor analysis model and more general models for Gaussian random
variables.  In this final section, we highlight two such problems: maximum
likelihood (ML) estimation and the study of singularities. Another
important direction is the design of test statistics from higher invariants
(pentads, septads, etc.) by computing moments of the Wishart distribution
(cf.~Section \ref{sec:tetradsstats}).

Since much statistical inference in factor analysis is currently based on
ML estimation, it seems natural to ask what computational algebra has to
say about computing ML estimators for factor analysis.  As a starting
point, we computed the \emph{maximum likelihood degree} for the one-factor
model with four observed variables.  The ML degree \citep{catanese:05} of a
statistical model is the number of nontrivial complex zeros of the critical
equations for generic data. We found that the factor analysis model
$\mathbf{F}_{4,1}$ has ML degree $57$.
In what follows, we describe how we discovered 
the number  $57$ to be the number of
complex zeros of the critical equations for ML
estimation in this model.

 Consider a matrix $\Psi \in F_{p,m}$ with decomposition
$\Psi=\Sigma+\Gamma$, where $\Sigma$ is diagonal with positive entries and
$\Gamma$ is positive semidefinite with $\rank(\Gamma)\le m$. 
Then we have
\[
\Psi^{-1}\quad = \quad \Sigma^{-1}-\Psi^{-1}\Gamma\Sigma^{-1}.
\]
The matrix $\Psi^{-1}\Gamma\Sigma^{-1}$ is symmetric of rank $\le m$.
Moreover, the fact that $\,\Psi-\Sigma\,=\,\Gamma \,$ is positive semidefinite
implies that $\,\Sigma^{-1}-\Psi^{-1}\,=\,\Psi^{-1}\Gamma\Sigma^{-1}\,$ is also
positive semidefinite.  Hence, the inverse of $\Psi \in F_{p,m}$ can be
written as
\[
\Psi^{-1} \quad = \quad T-KK^t,
\]
where $T$ is diagonal with positive entries and $K\in\RRR^{p\times m}$.
In the case of $p=4$ and $m=1$, this amounts to writing $\Psi^{-1}$ in
terms of $\tau=(\tau_1,\dots,\tau_4)$ and
$\kappa=(\kappa_1,\dots,\kappa_4)$ as
\begin{equation}
\label{eq:p4m1para}
\Psi^{-1}(\tau,\kappa) \quad = \quad
{\rm diag}(\tau) - \kappa \kappa^T.
\end{equation}

As in Section \ref{sec:tetradsstats}, let $\bar X$ and $S$ be the sample
mean vector and the sample covariance matrix computed from a sample of $N$
random vectors.  For ML estimation it is more convenient to work with the
matrix $\tilde S=(N-1)/N\cdot S$ instead of $S$.  The model of multivariate
normal distributions $\ND(\mu,\Psi)$ has the log-likelihood function
\begin{equation}
  \label{eq:fullloglik}
  \ell(\mu,\Psi) = -\frac{N}{2}\log \det(\Psi)
  -\frac{N}{2}\trace(\tilde S\Psi^{-1}) - \frac{N}{2}(\bar X-\mu)^t\Psi^{-1}(\bar X-\mu);
\end{equation}
compare \citet[Sect.~4.1.1]{mardia:79}.  This function is maximized in
$\mu$ by setting $\mu=\bar X$, which leads to the vanishing of the
quadratic form appearing as third term in (\ref{eq:fullloglik}).
Therefore, we can find the maximizer in $\Psi$ by maximizing the expression
\begin{equation}
  \label{eq:loglik}
  \log \det(\Psi^{-1})
  \,-\, \trace(\tilde S\Psi^{-1}).
\end{equation}
Here $\Psi $ runs over $F_{4,1}$.  
By plugging (\ref{eq:p4m1para}) into (\ref{eq:loglik}),
we can write this expression as a
 function of the eight unknowns
$\tau_1,\tau_2,\tau_3,\tau_4,
\kappa_1,\kappa_2,\kappa_3,\kappa_4 $.
Taking partial derivatives and setting them to zero,
we obtain a system of eight equations in eight unknowns.
These are the {\em likelihood equations} of the
factor analysis model $\mathbf{F}_{4,1}$ in rational function form:
\[
\begin{split}
\frac{1}{\det(\Psi^{-1}(\tau,\kappa))}\cdot\frac{
  \partial \det(\Psi^{-1}(\tau,\kappa))}{
  \partial \tau_i} &=\trace\bigg[\tilde S\frac{
  \partial \Psi^{-1}(\tau,\kappa)}{
  \partial \tau_i}\bigg], \qquad i=1,\dots,4,\\
\frac{1}{\det(\Psi^{-1}(\tau,\kappa))}\cdot\frac{
  \partial \det(\Psi^{-1}(\tau,\kappa))}{
  \partial \kappa_i} &=\trace\bigg[\tilde S\frac{
  \partial \Psi^{-1}(\tau,\kappa)}{
  \partial \kappa_i}\bigg], \qquad i=1,\dots,4.
\end{split}
\]
These equations can be made polynomial by multiplying through by
$\det(\Psi^{-1}(\tau,\kappa))$.  Clearing the denominator introduces many
additional solutions to the system, namely noninvertible matrices of the
form $\Psi^{-1}(\tau,\kappa)$.  However, these extraneous solutions can be
removed using an operation called {\em saturation\/}.  After saturation, we
discover that the solution set of the polynomial likelihood equations
consists of $57$ isolated (complex) points.  Clearly, these 57 solutions
come in pairs $(\tau,\pm\kappa)$; one solution has $\kappa=0$.  Further
work needs to be done to determine how many of these can be statistically
meaningful local maxima and to extend these results to larger models.

\smallskip

The second problem we wish to illustrate is that of singularities.  The
singularities of statistical models play an important, though
under-appreciated, role.  Models with singularities do not form curved
exponential families, which invalidates, for example, the theoretical basis
of model selection using information criteria like BIC \citep{geiger:01}.
Near a singularity, such criteria require correction terms
\citep{watanabe:01}.  This is especially important in factor analysis
because, among other singularities, each parameter space $F_{p,m}$ is
singular along $F_{p,m-1}$ making the selection of the number of factors
difficult.  The fact that $F_{p,m-1}$ is contained in the singular locus of
$F_{p,m}$ can be proved by observing that $V(I_{p,m})$ is the $m$-th secant
variety of $V(I_{p,1})$ and by appealing to general results in algebraic
geometry about singularities of secant varieties.

As a first step towards a better understanding of singularities, we
computed the singular loci of some of the small factor analysis models.
The computation of the singular locus is done, in {\tt Macaulay 2} or {\tt
  Singular}, by augmenting $I_{p,m}$ by the $c \times c$-minors of the
Jacobian matrix of any generating set of $I_{p,m}$, where $c = {\rm
  codim}(I_{p,m})$.

\begin{example}\label{ex:singone}
  For $p = 4$ and $m =1$, the singular locus consists of all matrices
  $\Psi$ which have at most one non-zero off-diagonal entry.  Thus, the
  singular locus consists of one symmetry class of covariance matrices.
  For instance, one type of matrix in the singular locus has the form
$$\begin{pmatrix}
\psi_{11} & \psi_{12} & 0 & 0 \\
\psi_{12} & \psi_{22} & 0 &  0\\
0 & 0 & \psi_{33} & 0 \\
0 & 0 & 0 & \psi_{44}
\end{pmatrix}.
$$
\end{example}

Example \ref{ex:singone} generalizes to an arbitrary number of observed
random variables when the number of factors is fixed at $m=1$.

\begin{proposition}
  A matrix $\Psi \in V_\RRR(I_{p,1})$ is a singularity of the one factor
  model if and only if $\Psi$ has at most one non-zero off diagonal entry.
\end{proposition} 

\begin{proof}
  This can be seen by computing derivatives of the tetrads and by using
  results in toric geometry.  We omit the details.
\end{proof}

For more factors, the situation is considerably more complicated.

\begin{example}
  Let $p = 5$ and $m = 2$.  In this case $F_{5,2}$ is the pentad
  hypersurface.  The singular locus of this hypersurface has dimension 11,
  and consists of two symmetry classes of singularities.  A representative
  from the first symmetry class is the set of matrices of the form
  $$
  \begin{pmatrix}
    \psi_{11} & 0 & 0 & 0 & 0 \\
    0 & \psi_{22} & \psi_{23} & \psi_{24} & \psi_{25} \\
    0 & \psi_{23} & \psi_{33} & \psi_{34} & \psi_{35} \\
    0 & \psi_{24} & \psi_{34} & \psi_{44} & \psi_{45} \\
    0 & \psi_{25} & \psi_{35} & \psi_{45} & \psi_{55} \\ 
  \end{pmatrix}.
  $$
  The second symmetry class is more complicated.  A representative set
  consists of those matrices $\Psi$ that satisfy all tetrads not involving
  $\psi_{12}$.  Note that this set of singular points contains $F_{5,1}$.
  In total there are five elements of the first symmetry class and ten
  elements in the second.  To apply algebraic geometry techniques to these
  singularities for model selection requires a careful analysis of the way
  the various singular sets intersect.
\end{example}

It is an open problem to determine the singular locus of the
factor analysis models $F_{p,m}$ in general. Even the dimension of that
singular locus is  unknown to us.

\medskip

Finally we stress that while the algebraic statistical study in this paper was confined
to factor analysis models, problems analogous to the ones described here
appear in other classes of Gaussian graphical models.  The study of such other
models, which need not involve hidden variables, opens up a broad range of
directions for future research.

\bigskip
\bigskip

\noindent {\bf Acknowledgements.}
Much of this work was conducted while Mathias Drton held a postdoctoral
position sponsored by the Center of Pure and Applied Mathematics at UC
Berkeley.  Mathias Drton also acknowledges support from the National
Science Foundation (DMS-0505612).  Seth Sullivant was supported by an NSF
Graduate Research Fellowship.  Bernd Sturmfels was also supported by the
NSF (DMS-0456960).

\bigskip
\bigskip

\small
\bibliographystyle{plainnat}

\bigskip
\bigskip
\bigskip
\bigskip
\bigskip

\noindent {\bf Authors' addresses:}

\bigskip

\noindent Mathias Drton, Department of Statistics,
University of Chicago,  Chicago, IL 60637, USA,
{\tt drton@galton.uchicago.edu}

\medskip

\noindent Bernd Sturmfels, Department of Mathematics,
University of California,  Berkeley, CA 94720, USA,
{\tt bernd@math.berkeley.edu}

\medskip

\noindent Seth Sullivant, Society of Fellows and Department of Mathematics,
Harvard University, Cambridge, MA 02138, USA
{\tt seths@math.harvard.edu}

\end{document}